\documentclass[a4paper,11pt]{article}
\usepackage{graphicx}
\usepackage{amsmath}
\usepackage{amsfonts}
\usepackage{amssymb}
\newtheorem{theorem}{Theorem}

\newtheorem{definition}[theorem]{Definition}
\newtheorem{example}[theorem]{Example}

\newtheorem{lemma}[theorem]{Lemma}
\newtheorem{notation}[theorem]{Notation}

\newtheorem{proposition}[theorem]{Proposition}
\newtheorem{remark}[theorem]{Remark}

\newenvironment{proof}[1][Proof]{\textbf{#1.} }{\ \rule{0.5em}{0.5em}}
\begin{document}

\title{Generalized Integral Operators and Schwartz Kernel Theorem}
\author{A. Delcroix\\Equipe Analyse Alg\'{e}brique Non Lin\'{e}aire\\\textit{Laboratoire Analyse, Optimisation, Contr\^{o}le}\\Facult\'{e} des sciences - Universit\'{e} des Antilles et de la Guyane\\97159\ Pointe-\`{a}-Pitre Cedex Guadeloupe}
\maketitle
\begin{abstract}
In connection with the classical Schwartz kernel theorem, we show that in the
framework of Colombeau generalized functions a large class of \ linear
mappings admit integral kernels. To do this, we need to introduce news spaces
of generalized functions with slow growth and the corresponding adapted linear
mappings. Finally, we show that in some sense Schwartz' result is contained in
our main theorem.
\end{abstract}
\tableofcontents

\noindent~\bigskip

\noindent\textbf{Mathematics Subject Classification (2000): 45P05, 46F05,
46F30, 47G10\smallskip}

\noindent\textbf{Keywords:} Schwartz kernel theorem, generalized functions,
integral operator, distributions.\bigskip

\section{Introduction}

It is well known that the framework of Schwartz distributions is not suitable
for setting and solving many differential or integral problems with singular
coefficients or data. A natural approach to overcome this difficulty consists
in replacing the given problem by a one-parameter family of smooth
problems.\ This is done in most theories of generalized functions and, for
example, in Colombeau simplified theory which we are going to use in the
sequel. (For details, see the monographies \cite{Col1}, \cite{GKOS},
\cite{Ober1} and the references therein.)

In this paper, we continue the investigations in the field of generalized
integral operators initiated by the pioneering work of D.\ Scarpalezos
\cite{Scarpa1}, and carried on by J.-F.\ Colombeau (personal communications
and \cite{BCD}) in view of applications to physics and by\ C.\ Garetto
\emph{et alii} (\cite{GaGrOb}) with applications to pseudo differential
operators theory and questions of regularity.\medskip

More precisely, the following results holds: Every $H$ belonging to
$\mathcal{G}\left(  \mathbb{R}^{m}\times\mathbb{R}^{n}\right)  $ defines a
linear operator from $\mathcal{G}_{C}\left(  \mathbb{R}^{n}\right)  $ to
$\mathcal{G}\left(  \mathbb{R}^{m}\right)  $ by the formula
\[
\widetilde{H}:\mathcal{G}_{C}\left(  \mathbb{R}^{n}\right)  \rightarrow
\mathcal{G}\left(  \mathbb{R}^{m}\right)  ,\;\;\;\;f\mapsto\widetilde
{H}(f)\text{ with }\widetilde{H}(f)(x)=\left[  \int H_{\varepsilon
}(x,y)f_{\varepsilon}(y)\,\mathrm{d}y\right]  ,
\]
where $\left(  H_{\varepsilon}\right)  _{\varepsilon}$ (\textit{resp.}
$\left(  f_{\varepsilon}\right)  _{\varepsilon}$) is any representative of $H$
(\textit{resp.} $f$) and \ $\left[  \,\cdot\,\right]  $ is the class of an
element in $\mathcal{G}\left(  \mathbb{R}^{d}\right)  $. ($\mathcal{G}\left(
\mathbb{R}^{d}\right)  $ denotes the usual quotient space of Colombeau
simplified generalized functions, while $\mathcal{G}_{C}\left(  \mathbb{R}%
^{d}\right)  $ is the subspace of elements of $\mathcal{G}\left(
\mathbb{R}^{d}\right)  $ compactly supported: See section \ref{GSTSecFW} for
the mathematical framework.)

Conversely, in the distributional case, the well known Schwartz kernel theorem
asserts that each linear map $\Lambda$ from $\mathcal{D}\left(  \mathbb{R}%
^{n}\right)  $ to $\mathcal{D}^{\prime}\left(  \mathbb{R}^{m}\right)  $
continuous for the strong topology of $\mathcal{D}^{\prime}$ can be
represented by a kernel $K\in\mathcal{D}^{\prime}\left(  \mathbb{R}^{m}%
\times\mathbb{R}^{n}\right)  $ that is
\[
\forall f\in\mathcal{D}\left(  \mathbb{R}^{n}\right)  ,\;\;\forall\varphi
\in\mathcal{D}\left(  \mathbb{R}^{m}\right)  ,\;\;\left(  \Lambda\left(
f\right)  ,\varphi\right)  =\left(  K,\varphi\otimes f\right)  .
\]

Let us recall here that $\mathcal{D}\left(  \mathbb{R}^{n}\right)  $ is
embedded in $\mathcal{G}_{C}\left(  \mathbb{R}^{n}\right)  $ and
$\mathcal{D}^{\prime}\left(  \mathbb{R}^{m}\right)  $ in $\mathcal{G}\left(
\mathbb{R}^{m}\right)  $: In the spirit of Schwartz theorem, we prove that in
the framework of Colombeau generalized functions any net of linear maps
$\left(  L_{\varepsilon}:\mathcal{D}\left(  \mathbb{R}^{n}\right)
\rightarrow\mathrm{C}^{\infty}\left(  \mathbb{R}^{m}\right)  \right)
_{\varepsilon}$ satisfying some growth property with respect to the parameter
$\varepsilon$ (the \emph{strongly moderate nets}) gives rise to a linear map
$L:\mathcal{G}_{C}\left(  \mathbb{R}^{n}\right)  \rightarrow\mathcal{G}\left(
\mathbb{R}^{m}\right)  $ which can be represented as an integral operator.
This means that there exists a generalized function $H_{L}\in\mathcal{G}%
\left(  \mathbb{R}^{m}\times\mathbb{R}^{n}\right)  $ such that
\[
L(f)=\int H_{L}(\cdot,y)f(y)\,\mathrm{d}y
\]
for any $f$ belonging to a convenient subspace of $\mathcal{G}_{C}\left(
\mathbb{R}^{n}\right)  $.

Moreover, this result is strongly related to Schwartz Kernel theorem in the
following sense.\ We can associate to each linear operator $\Lambda
:\mathcal{D}\left(  \mathbb{R}^{n}\right)  \rightarrow\mathcal{D}^{\prime
}\left(  \mathbb{R}^{m}\right)  $ satisfying the hypothesis above mentioned\ a
strongly moderate map $L_{\Lambda}$ and consequently a kernel $H_{L_{\Lambda}%
}\in\mathcal{G}\left(  \mathbb{R}^{m}\times\mathbb{R}^{n}\right)  $ with the
following equality property: For all $f$ in $\mathcal{D}\left(  \mathbb{R}%
^{n}\right)  $, $\Lambda\left(  f\right)  $ and $\widetilde{H}_{L_{\Lambda}%
}\left(  f\right)  $ are equal in the generalized distribution sense
\cite{NePiSc} that is, for all $k\in\mathbb{N}$ and $\left(  H_{L_{\Lambda
,\varepsilon}}\right)  _{\varepsilon}$ representative of $H_{L_{\Lambda}}$,
\[
\forall\varphi\in\mathcal{D}\left(  \mathbb{R}^{m}\right)  ,\;\;\left\langle
\Lambda\left(  f\right)  ,\Phi\right\rangle -%
{\textstyle\int}
\left(
{\textstyle\int}
H_{L_{\Lambda,\varepsilon}}\left(  x,y\right)  f\left(  y\right)
\,\mathrm{d}y\right)  \varphi\left(  x\right)  \,\mathrm{d}x=\mathrm{O}\left(
\varepsilon^{k}\right)  \text{, for }\varepsilon\rightarrow0.
\]

The paper can be divided in two parts.\ The first part, formed by section
\ref{GSTSecFW} and section \ref{GSTSSPrelST}, introduces all the material
which is needed in the sequel.\ We mention here in particular the notion of
\textit{spaces of generalized functions with slow growth}, which are subspaces
of the usual space $\mathcal{G}\left(  \mathbb{R}^{d}\right)  $ with
additional limited growth property with respect to the parameter $\varepsilon
$. Lemma \ref{LmnMTh4} shows one feature of those spaces (used for the proof
of the main results): Convolution admits on them as unit some special $\delta
$-nets, whereas with result is false in $\mathcal{G}\left(  \mathbb{R}%
^{d}\right)  $. The second part, formed by the two last sections, is devoted
to the definition of strongly moderate nets, the statement of the main results
and their proofs.

\section{Colombeau type algebras\label{GSTSecFW}}

\subsection{The sheaf of Colombeau simplified algebras}

Let $\mathrm{C}^{\infty}$ be the sheaf of complex valued smooth functions on
$\mathbb{R}^{d}$ ($d\in\mathbb{N}$) with the usual topology of uniform
convergence.\ For every open set $\Omega$ of $\mathbb{R}^{d}$, this topology
can be described by the family of semi norms
\[
p_{K,l}(f)=\sup_{\left|  \alpha\right|  \leq l,K\Subset\Omega}\left|
\partial^{\alpha}f\left(  x\right)  \right|
\]
where the notation $K\Subset\Omega$ means that\ the set $K$ is a compact set
included in $\Omega$.

Let us set\smallskip

\noindent$\mathcal{F}\left(  \mathrm{C}^{\infty}\left(  \Omega\right)
\right)  =\left\{  \left(  f_{\varepsilon}\right)  _{\varepsilon}\in
\mathrm{C}^{\infty}\left(  \Omega\right)  ^{\left(  0,1\right]  }\,\left|
\,\forall l\in\mathbb{N},\;\forall K\Subset\Omega,\;\exists q\in
\mathbb{N},\;p_{K,l}\left(  f_{\varepsilon}\right)  =\mathrm{O}\left(
\varepsilon^{-q}\right)  \;\mathrm{for}\;\varepsilon\rightarrow0\right.
\right\}  ,$\smallskip

\noindent$\mathcal{N}\left(  \mathrm{C}^{\infty}\left(  \Omega\right)
\right)  =\left\{  \left(  f_{\varepsilon}\right)  _{\varepsilon}\in
\mathrm{C}^{\infty}\left(  \Omega\right)  ^{\left(  0,1\right]  }\,\left|
\,\forall l\in\mathbb{N},\;\forall K\Subset\Omega,\;\forall p\in
\mathbb{N},\;p_{K,l}\left(  f_{\varepsilon}\right)  =\mathrm{O}\left(
\varepsilon^{p}\right)  \;\mathrm{for}\;\varepsilon\rightarrow0\right.
\right\}  .$\smallskip

\begin{lemma}
\cite{JAM1} and \cite{JAM2}\newline $i$.~The functor $\mathcal{F}%
:\Omega\rightarrow\mathcal{F}\left(  \mathrm{C}^{\infty}\left(  \Omega\right)
\right)  $ defines a sheaf of subalgebras of the sheaf $\left(  \mathrm{C}%
^{\infty}\right)  ^{\left(  0,1\right]  }$\newline $ii$.~The functor
$\mathcal{N}:\Omega\rightarrow\mathcal{N}\left(  \mathrm{C}^{\infty}\left(
\Omega\right)  \right)  $ defines a sheaf of ideals of the sheaf $\mathcal{F}$.
\end{lemma}

We shall note prove in detail this lemma but quote the two mains arguments:

\noindent$i$.~For each open subset $\Omega$ of $X$, the family of seminorms
$\left(  p_{K,l}\right)  $ related to $\Omega$ is compatible with the
algebraic structure of $\mathcal{E}\left(  \Omega\right)  ;$\ In particular:
\[
\forall l\in\mathbb{N},\;\forall K\Subset\Omega,\;\;\exists C\in\mathbb{R}%
_{+}^{\ast},\;\;\forall\left(  f,g\right)  \in\left(  \mathrm{C}^{\infty
}\left(  \Omega\right)  \right)  ^{2}\;\;p_{K,l}\left(  fg\right)  \leq
Cp_{K,l}\left(  f\right)  p_{K,l}\left(  g\right)  ,
\]

\noindent$ii$.~For two open subsets $\Omega_{1}\subset\Omega_{2}$ of
$\mathbb{R}^{d}$, the family of seminorms $\left(  p_{K,l}\right)  $ related
to $\Omega_{1}$ is included in the family of seminorms related to $\Omega_{2}$
and
\[
\forall l\in\mathbb{N},\;\forall K\Subset\Omega_{1},\;\;\forall f\in
\mathrm{C}^{\infty}\left(  \Omega_{2}\right)  ,\;\;p_{K,l}\left(  f_{\left|
\Omega_{1}\right.  }\right)  =p_{K,l}\left(  f\right)  .
\]

\begin{definition}
The sheaf of factor algebras
\[
\mathcal{G}=\mathcal{F}\left(  \mathrm{C}^{\infty}\left(  \cdot\right)
\right)  /\mathcal{N}\left(  \mathrm{C}^{\infty}\left(  \cdot\right)
\right)
\]
is called the sheaf of \emph{Colombeau type algebras}.
\end{definition}

The sheaf $\mathcal{G}$ turns to be a sheaf of differential algebras and a
sheaf of modulus on the factor ring $\overline{\mathbb{C}}=\mathcal{F}\left(
\mathbb{C}\right)  /\mathcal{N}\left(  \mathbb{C}\right)  $ with
\begin{align*}
\mathcal{F}\left(  \mathbb{K}\right)   &  =\left\{  \left(  r_{\varepsilon
}\right)  _{\varepsilon}\in\mathbb{K}^{\left(  0,1\right]  }\,\left|
\,\exists q\in\mathbb{N},\;\left|  r_{\varepsilon}\right|  =\mathrm{O}\left(
\varepsilon^{-q}\right)  \;\mathrm{for}\;\varepsilon\rightarrow0\right.
\right\}  ,\\
\mathcal{N}\left(  \mathbb{K}\right)   &  =\left\{  \left(  r_{\varepsilon
}\right)  _{\varepsilon}\in\mathbb{K}^{\left(  0,1\right]  }\,\left|
\,\forall p\in\mathbb{N},\;\left|  r_{\varepsilon}\right|  =\mathrm{O}\left(
\varepsilon^{p}\right)  \;\mathrm{for}\;\varepsilon\rightarrow0\right.
\right\}  ,
\end{align*}
with $\mathbb{K}=\mathbb{C}$ or $\mathbb{K}=\mathbb{R},\;\mathbb{R}_{+}$.

\begin{notation}
In the sequel we shall note, as usual, $\mathcal{G}\left(  \Omega\right)  $
instead of $\mathcal{G}\left(  \mathrm{C}^{\infty}\left(  \Omega\right)
\right)  $ the algebra of generalized functions on $\Omega$. For $\left(
f_{\varepsilon}\right)  _{\varepsilon}\in\mathcal{F}\left(  \mathrm{C}%
^{\infty}\left(  \Omega\right)  \right)  $, $\left[  \left(  f_{\varepsilon
}\right)  _{\varepsilon}\right]  $ will be its class in $\mathcal{G}\left(
\Omega\right)  $.
\end{notation}

\subsection{Generalized functions with compact supports\label{GSTSSCCompSupp}}

Let us mention here some remarks about generalized functions with compact
supports, which will be useful in the sequel.

As $\mathcal{G}$ is a sheaf, the notion of support of a section $f\in
\mathcal{G}\left(  \Omega\right)  $ make sense. The following definition will
be sufficient for this paper.

\begin{definition}
The support of a generalized function $f\in\mathcal{G}\left(  \Omega\right)  $
is the complement in $\Omega$ of the largest open subset of $\Omega$ where $f
$ is null.
\end{definition}

\begin{notation}
We denote by $\mathcal{G}_{C}\left(  \Omega\right)  $ the subset of
$\mathcal{G}\left(  \Omega\right)  $ of elements with compact supports.
\end{notation}

\begin{lemma}
\label{LmnSuppK}Every $f\in\mathcal{G}_{C}$ has a representative $\left(
f_{\varepsilon}\right)  _{\varepsilon}$, such that each $f_{\varepsilon}$ has
the same compact support.
\end{lemma}

There is an other way to introduce generalized functions with compact support
more natural in the sequel. We start from the algebra $\mathcal{D}\left(
\Omega\right)  $ considered as the inductive limit of
\[
\mathcal{D}_{j}\left(  \Omega\right)  =\mathcal{D}_{K_{j}}\left(
\Omega\right)  =\left\{  f\in\mathcal{D}\left(  \Omega\right)  \,\left|
\,\operatorname*{supp}f\subset K_{j}\right.  \right\}
\]
where:

\begin{description}
\item [i.]$\left(  K_{j}\right)  _{j\in\mathbb{N}}$ is an increasing sequence
of relatively compact subsets exhausting $\Omega$, with $K_{j}\subset
\overset{\circ}{K}_{j+1};$

\item[ii.] $\mathcal{D}_{j}\left(  \Omega\right)  $ is endowed with the family
of semi norms
\[
p_{j,l}\left(  f\right)  =\sup_{\left|  \alpha\right|  \leq l,x\in K_{j}%
}\left|  \partial^{\alpha}f\left(  x\right)  \right|  .
\]
\end{description}

The topology on $\mathcal{D}\left(  \Omega\right)  $ does not depend on the
particular choice of the sequence $\left(  K_{j}\right)  _{j\in\mathbb{N}}$.
Construction of spaces of generalized functions based on projective or
inductive limits have already been considered (see e.g. \cite{DHPV1},
\cite{PScamb}). We just recall it briefly here. Let fix $\left(  K_{j}\right)
_{j\in\mathbb{N}}$ a sequence of compact sets satisfying \textbf{i.} and set
\begin{gather}
\mathcal{F}\left(  \mathcal{D}\left(  \Omega\right)  \right)  =\cup
_{j\in\mathbb{N}}\mathcal{F}_{j}\left(  \Omega\right)  \text{ }%
\label{ModCompSup}\\
\text{with\ }\mathcal{F}_{j}\left(  \Omega\right)  =\left\{  \left(
f_{\varepsilon}\right)  _{\varepsilon}\in\mathcal{D}_{j}\left(  \Omega\right)
^{\left(  0,1\right]  }\,\left|  \,\forall l\in\mathbb{N},\;\exists
q\in\mathbb{N},\;p_{j,l}\left(  f_{\varepsilon}\right)  =\mathrm{O}\left(
\varepsilon^{-q}\right)  \;\mathrm{for}\;\varepsilon\rightarrow0\right.
\right\} \nonumber\\
\mathcal{N}\left(  \mathcal{D}\left(  \Omega\right)  \right)  =\cup
_{n\in\mathbb{N}}\mathcal{N}_{n}\left(  \Omega\right)  \text{ ,}\nonumber\\
\text{with\ }\mathcal{N}_{j}\left(  \Omega\right)  =\left\{  \left(
f_{\varepsilon}\right)  _{\varepsilon}\in\mathcal{D}_{j}\left(  \Omega\right)
^{\left(  0,1\right]  }\,\left|  \,\forall l\in\mathbb{N},\;\forall
p\in\mathbb{N},\;p_{j,l}\left(  f_{\varepsilon}\right)  =\mathrm{O}\left(
\varepsilon^{p}\right)  \;\mathrm{for}\;\varepsilon\rightarrow0\right.
\right\}  .\nonumber
\end{gather}

With those definitions, we have:

\begin{lemma}
$\mathcal{F}\left(  \mathcal{D}\left(  \Omega\right)  \right)  $ is a
subalgebra of $\mathcal{D}\left(  \Omega\right)  ^{\left(  0,1\right]  }$ and
$\mathcal{N}\left(  \mathcal{D}\left(  \Omega\right)  \right)  $ an ideal of
$\mathcal{F}\left(  \mathcal{D}\left(  \Omega\right)  \right)  $.
\end{lemma}

The factor space $\mathcal{G}_{\mathcal{D}}\left(  \Omega\right)
=\mathcal{F}\left(  \mathcal{D}\left(  \Omega\right)  \right)  /\mathcal{N}%
\left(  \mathcal{D}\left(  \Omega\right)  \right)  $ appears to be a natural
space of generalized functions with compact supports. The algebra
$\mathcal{G}_{\mathcal{D}}\left(  \Omega\right)  $ does not depend on the
particular choice of the sequence $\left(  K_{j}\right)  _{j\in\mathbb{N}}%
$\ Moreover, due to the properties of the family $\left(  p_{j,l}\right)  $ we have:

\begin{lemma}
\label{LmnKDKC1}The spaces $\mathcal{G}_{\mathcal{D}}\left(  \Omega\right)  $
and $\mathcal{G}_{C}\left(  \Omega\right)  $ are isomorphic.
\end{lemma}

\begin{proof}
The fundamental property involved is the following: for all $j\in\mathbb{N}$
and all $\left(  f_{\varepsilon}\right)  _{\varepsilon}\in\mathcal{F}%
_{j}\left(  \Omega\right)  $ we have
\begin{equation}
\forall l\in\mathbb{N},\;\;\forall j^{\prime}\leq j,\;\;\forall j^{\prime
\prime}\geq j,\;\;\;p_{j^{\prime},l}\left(  f_{\varepsilon}\right)  \leq
p_{_{j,l}}\left(  f_{\varepsilon}\right)  =p_{j^{\prime\prime},l}\left(
f_{\varepsilon}\right)  . \label{LmnKDKC2}%
\end{equation}
The last equality is true since $\operatorname*{supp}f\subset K_{j}\subset
K_{j^{\prime\prime}}$, for all $j^{\prime\prime}\geq j$.

Relation (\ref{LmnKDKC2}) implies that $\mathcal{F}\left(  \mathcal{D}\left(
\Omega\right)  \right)  \subset\mathcal{F}\left(  \mathrm{C}^{\infty}\left(
\mathbb{\Omega}\right)  \right)  $ and $\mathcal{N}\left(  \mathcal{D}\left(
\Omega\right)  \right)  \subset\mathcal{N}\left(  \mathrm{C}^{\infty}\left(
\mathbb{\Omega}\right)  \right)  $. Let us show the first inclusion.\ Consider
$\left(  f_{\varepsilon}\right)  _{\varepsilon}$ in some $\mathcal{F}%
_{j}\left(  \Omega\right)  $. Then, for all $l\in\mathbb{N}$, there exists
$q\in\mathbb{N}$ such that: $\forall j\in\mathbb{N}$, $p_{_{j,l}}\left(
f_{\varepsilon}\right)  =\mathrm{O}\left(  \varepsilon^{-q}\right)  $ for
$\varepsilon\rightarrow0$. It follows that $\forall K\Subset\Omega$,
$p_{_{K,l}}\left(  f_{\varepsilon}\right)  =\mathrm{O}\left(  \varepsilon
^{-q}\right)  $ since the sequence $\left(  K_{j}\right)  _{j\in\mathbb{N}}$
exhausts $K$.

Those two inclusions implies that the map
\[
\iota:\;\;\mathcal{G}_{\mathcal{D}}\left(  \Omega\right)  \rightarrow
\mathcal{G}\left(  \Omega\right)  ,\;\;\;\;\;\;\left(  f_{\varepsilon}\right)
_{\varepsilon}+\mathcal{N}\left(  \mathcal{D}\left(  \Omega\right)  \right)
\mapsto\left(  f_{\varepsilon}\right)  _{\varepsilon}+\mathcal{N}\left(
\mathrm{C}^{\infty}\left(  \mathbb{\Omega}\right)  \right)
\]
is well defined with $\iota\left(  \mathcal{G}_{\mathcal{D}}\left(
\Omega\right)  \right)  \subset\mathcal{G}_{C}\left(  \Omega\right)  $.

It remains to show that the map $\iota$ is bijective.\ Indeed, if $\left(
f_{\varepsilon}\right)  _{\varepsilon}\in\mathcal{N}\left(  \mathrm{C}%
^{\infty}\left(  \mathbb{\Omega}\right)  \right)  $ with $\left(
f_{\varepsilon}\right)  _{\varepsilon}\in\mathcal{F}_{j}\left(  \Omega\right)
$, we have $\left(  f_{\varepsilon}\right)  _{\varepsilon}\in\mathcal{N}%
_{j}\left(  \Omega\right)  $ and $\left(  f_{\varepsilon}\right)
_{\varepsilon}\in\mathcal{N}\left(  \mathcal{D}\left(  \Omega\right)  \right)
$. Injectivity follows. Conversely, take $g\in\mathcal{G}_{C}\left(
\Omega\right)  $.\ According to lemma \ref{LmnSuppK}, there exists a compact
$K$ and a representative $\left(  g_{\varepsilon}\right)  _{\varepsilon}$ of
$g$ such that $\operatorname*{supp}g_{\varepsilon}\subset K$, for all
$\varepsilon$. We observe that $K$ is included in some $K_{j}$ and then that
$\left(  g_{\varepsilon}\right)  _{\varepsilon}\in\mathcal{F}_{j}\left(
\Omega\right)  $.\ finally, $\iota\left(  \left(  g_{\varepsilon}\right)
_{\varepsilon}+\mathcal{N}\left(  \mathcal{D}\left(  \Omega\right)  \right)
\right)  =g$.
\end{proof}

\subsection{Embeddings}

The space $\mathrm{C}^{\infty}\left(  \mathbb{R}^{d}\right)  $ ($d\in
\mathbb{N}$) is embedded in $\mathcal{G}\left(  \mathbb{R}^{d}\right)  $ by
the canonical map
\[
\sigma:\;\mathrm{C}^{\infty}\left(  \mathbb{R}^{d}\right)  \rightarrow
\mathcal{G}\left(  \mathbb{R}^{d}\right)  \;\;\;f\rightarrow\left(
f_{\varepsilon}\right)  _{\varepsilon}+\mathcal{N}\left(  \mathrm{C}^{\infty
}\left(  \mathbb{R}^{d}\right)  \right)  \text{, \ \ with }f_{\varepsilon
}=f\text{ for all }\varepsilon\in\left(  0,1\right]
\]
which is an injective homomorphism of algebras.

Moreover, the construction of $\mathcal{G}\left(  \mathbb{R}^{d}\right)  $
permits to embed the space $\mathcal{D}^{\prime}\left(  \mathbb{R}^{d}\right)
$ by means of convolution with suitable mollifiers. We follow in this paper
the ideas of \cite{NePiSc}.

\begin{lemma}
\label{LmnGoodM}There exists a net of mollifiers $\left(  \theta_{\varepsilon
}\right)  _{\varepsilon}\in\mathcal{D}\left(  \mathbb{R}^{d}\right)  ^{\left(
0,1\right]  }$ for all $\varepsilon$, such that for all $k\in\mathbb{N}$%
\begin{gather}
\int\theta_{\varepsilon}\left(  x\right)  \,\mathrm{d}x=1+\mathrm{O}\left(
\varepsilon^{k}\right)  \;\text{for}\;\varepsilon\rightarrow
0,\label{LmnGoodM1}\\
\forall m\in\mathbb{N}^{d}\backslash\left\{  0\right\}  ,\;\int x^{m}%
\theta_{\varepsilon}\left(  x\right)  \,\mathrm{d}x=\mathrm{O}\left(
\varepsilon^{k}\right)  \;\text{for}\;\varepsilon\rightarrow0.
\label{LmnGoodM2}%
\end{gather}
\end{lemma}

Such a net is built in the following way: Consider $\rho\in\mathcal{S}\left(
\mathbb{R}^{d}\right)  $ such that $\int\rho\left(  x\right)  \,\mathrm{d}%
x=1$,\ $\int x^{m}\rho\left(  x\right)  \,\mathrm{d}x=0$ for all
$m\in\mathbb{N}^{d}\backslash\left\{  0\right\}  $ and $\kappa\in
\mathcal{D}\left(  \mathbb{R}^{d}\right)  $ such that $0\leq\kappa\leq1$,
$\kappa=1$ on $\left[  -1,1\right]  ^{d}$ and $\kappa=0$ on $\mathbb{R}%
^{d}\backslash\left[  -2,2\right]  ^{d}$.\ Then \ $\left(  \theta
_{\varepsilon}\right)  _{\varepsilon}$ defined by
\[
\forall\varepsilon\in\left(  0,1\right]  ,\;\;\forall x\in\mathbb{R}%
^{d},\;\;\;\;\theta_{\varepsilon}\left(  x\right)  =\frac{1}{\varepsilon^{d}%
}\rho\left(  \frac{x}{\varepsilon}\right)  \kappa\left(  x\left|
\ln\varepsilon\right|  \right)  \text{ }%
\]
satisfies conditions of lemma \ref{LmnGoodM}.

\begin{proposition}
\label{LmnEmbed}With notations of lemma \ref{LmnGoodM}, the map
\[
\iota:\;\mathcal{D}^{\prime}\left(  \mathbb{R}^{d}\right)  \rightarrow
\mathcal{G}\left(  \mathbb{R}^{d}\right)  \;\;\;T\mapsto\left(  T\ast
\theta_{\varepsilon}\right)  _{\varepsilon}+\mathcal{N}\left(  \mathrm{C}%
^{\infty}\left(  \mathbb{R}^{d}\right)  \right)
\]
is an injective homomorphism of vector spaces.\ Moreover $\iota_{\left|
\mathrm{C}^{\infty}\left(  \Omega\right)  \right.  }=\sigma$.
\end{proposition}

This proposition asserts that the following diagram is commutative:

\begin{center}
$%
\begin{array}
[c]{ccc}%
\mathrm{C}^{\infty}\left(  \mathbb{R}^{d}\right)  & \longrightarrow &
\mathcal{D}^{\prime}\left(  \mathbb{R}^{d}\right) \\
& \searrow\!\sigma & \downarrow\iota\\
&  & \mathcal{G}\left(  \mathbb{R}^{d}\right)
\end{array}
$
\end{center}

\subsection{Generalized Integral operators}

We collect here results about generalized integral operators. We refer the
reader to \cite{BCD} and \cite{GaGrOb} for details.

\begin{definition}
\label{DefGIOP}Let $H$ be in $\mathcal{G}(\mathbb{R}^{m}\times\mathbb{R}^{n}%
)$.The integral operator of Kernel $H$ is the map $\widetilde{H}$ defined by
\[
\widetilde{H}:\mathcal{G}_{C}\left(  \mathbb{R}^{n}\right)  \rightarrow
\mathcal{G}\left(  \mathbb{R}^{m}\right)  :f\mapsto\widetilde{H}\left(
f\right)  \text{ with }\widetilde{H}\left(  f\right)  =\left[  \left(
x\mapsto\int H_{\varepsilon}(x,y)f(y)\,\mathrm{d}y\right)  _{\varepsilon
}\right]
\]
where $\left(  H_{\varepsilon}\right)  _{\varepsilon}$ is any representative
of $H$.
\end{definition}

Note that in the above mentioned references, the generalized function $H$
satisfies some additive condition such as being properly supported.\ This
assumption is not needed in this paper, since we consider operators on
$\mathcal{G}_{C}\left(  \mathbb{R}^{n}\right)  $: the integral which appears
in definition\ \ref{DefGIOP} is performed on a compact set.

\begin{proposition}
\label{PropOpInt}With the notations of definition \ref{DefGIOP} the operator
$\widetilde{H}$ defines a linear mapping from $\mathcal{G}_{C}\left(
\mathbb{R}^{n}\right)  $ to $\mathcal{G}\left(  \mathbb{R}^{m}\right)  $
continuous for the respective sharp topologies of $\mathcal{G}_{C}\left(
\mathbb{R}^{n}\right)  $ and $\mathcal{G}\left(  \mathbb{R}^{m}\right)
$.\newline Moreover the map
\[
\mathcal{G}(\mathbb{R}^{m}\times\mathbb{R}^{n})\rightarrow\mathcal{L}\left(
\mathcal{G}_{C}\left(  \mathbb{R}^{n}\right)  ,\mathcal{G}\left(
\mathbb{R}^{m}\right)  \right)  \;\;\;\;\;H\mapsto\widetilde{H}%
\]
is injective.
\end{proposition}

In other words, the map $\widetilde{H}$ is characterized by the kernel $H$%
\[
\widetilde{H}=0\text{ in }\mathcal{L}\left(  \mathcal{G}_{C}\left(
\mathbb{R}^{n}\right)  ,\mathcal{G}\left(  \mathbb{R}^{m}\right)  \right)
\Leftrightarrow H=0\text{ in }\mathcal{G}(\mathbb{R}^{m}\times\mathbb{R}%
^{n}).
\]

\section{Spaces of generalized functions with slow growth\label{GSTSSPrelST}}

In the sequel, we need to consider some subspaces of $\mathcal{G}\left(
\Omega\right)  $ with restrictive conditions of growth with respect to
$1/\varepsilon$ when the $l$ index of the families of seminorms is involved,
that is the index related to derivatives. We show that these spaces give a
good framework for extension of linear maps and for convolution of generalized
functions. These are essential properties for our result.

\subsection{Definitions}

Set
\begin{align*}
\mathcal{F}_{\mathcal{L}_{0}}\left(  \mathrm{C}^{\infty}\left(  \Omega\right)
\right)   &  =\left\{  \left(  f_{\varepsilon}\right)  _{\varepsilon}%
\in\mathcal{F}\left(  \Omega\right)  ^{\left(  0,1\right]  }\,\left|  \forall
K\Subset\Omega,\;\exists q\in\mathbb{N}^{\mathbb{N}},\mathrm{\;with}\text{
}\lim_{l\rightarrow+\infty}\left(  q(l)/l\right)  =0\right.  \right. \\
&
\;\;\;\;\;\;\;\;\;\;\;\;\;\;\;\;\;\;\;\;\;\;\;\;\;\;\;\;\;\;\;\;\;\;\;\;\;\;\;\;\;\;\left.
\forall l\in\mathbb{N},\;\;\;p_{K,l}\left(  f_{\varepsilon}\right)
=\mathrm{O}\left(  \varepsilon^{-q(l)}\right)  \;\mathrm{for}\;\varepsilon
\rightarrow0\right\}  .\\
\mathcal{F}_{\mathcal{L}_{1}}\left(  \mathrm{C}^{\infty}\left(  \Omega\right)
\right)   &  =\left\{  \left(  f_{\varepsilon}\right)  _{\varepsilon}%
\in\mathcal{F}\left(  \Omega\right)  ^{\left(  0,1\right]  }\,\left|  \forall
K\Subset\Omega,\;\exists q\in\mathbb{N}^{\mathbb{N}},\mathrm{\;with}\text{
}\underset{l\rightarrow+\infty}{\lim\sup}\left(  q(l)/l\right)  <1\right.
\right. \\
&
\;\;\;\;\;\;\;\;\;\;\;\;\;\;\;\;\;\;\;\;\;\;\;\;\;\;\;\;\;\;\;\;\;\;\;\;\;\;\;\;\;\;\left.
\forall l\in\mathbb{N},\;\;\;p_{K,l}\left(  f_{\varepsilon}\right)
=\mathrm{O}\left(  \varepsilon^{-q(l)}\right)  \;\mathrm{for}\;\varepsilon
\rightarrow0\right\}  .
\end{align*}

\begin{lemma}
\label{SpLGr1}~\newline $i.~\mathcal{F}_{\mathcal{L}_{0}}\left(
\mathrm{C}^{\infty}\left(  \Omega\right)  \right)  $ is a subalgebra of
$\mathcal{F}\left(  \mathrm{C}^{\infty}\left(  \Omega\right)  \right)
$.$\newline $ii.$~\mathcal{F}_{\mathcal{L}_{1}}\left(  \mathrm{C}^{\infty
}\left(  \Omega\right)  \right)  $ is a submodulus of $\mathcal{F}\left(
\mathrm{C}^{\infty}\left(  \Omega\right)  \right)  .$
\end{lemma}

\begin{proof}
Take $\left(  f_{\varepsilon}\right)  _{\varepsilon}$ and $\left(
g_{\varepsilon}\right)  _{\varepsilon}$ in $\mathcal{F}_{\mathcal{L}_{0}%
}\left(  \mathrm{C}^{\infty}\left(  \Omega\right)  \right)  $ (\textit{resp.}
$\mathcal{F}_{\mathcal{L}_{1}}\left(  \mathrm{C}^{\infty}\left(
\Omega\right)  \right)  $), $K\Subset\Omega$, $q_{f}$ and $q_{g}$ the
corresponding sequences with $\lim_{l\rightarrow+\infty}\left(  q_{h}%
(l)/l\right)  =0$ (\textit{resp}. $r_{h}=\lim\!\sup_{l\rightarrow+\infty
}\left(  q_{h}(l)/l\right)  <1$) for $h=f,\,g$. Define
\[
q\left(  \cdot\right)  =\max\left(  q_{f}\left(  \cdot\right)  ,q_{g}\left(
\cdot\right)  \right)  \;\;\;\;\;\ r=\max\left(  r_{f},r_{g}\right)  <1\text{,
for the \textit{resp}. case.}%
\]
For $h=f,\,g$ we have $p_{K,l}\left(  h_{\varepsilon}\right)  =\mathrm{O}%
\left(  \varepsilon^{-q(l)}\right)  $ for\ $\varepsilon\rightarrow0$ and
$p_{K,l}\left(  f_{\varepsilon}+g_{\varepsilon}\right)  =\mathrm{O}\left(
\varepsilon^{-q(l)}\right)  $ for\ $\varepsilon\rightarrow0$.

For $\left(  c_{\varepsilon}\right)  _{\varepsilon}\in\mathcal{F}\left(
\mathbb{R}\right)  $, there exists $q_{c}$ such that $\left|  c_{\varepsilon
}\right|  =\mathrm{O}\left(  \varepsilon^{-q_{c}}\right)  $. Then
$p_{K,l}\left(  c_{\varepsilon}f_{\varepsilon}\right)  =\mathrm{O}\left(
\varepsilon^{-\left(  q_{c}+q(l)\right)  }\right)  $ with $\lim_{l\rightarrow
+\infty}\left(  \left(  q_{c}+q(l)\right)  /l\right)  =0$ (\textit{resp.}
$\lim\!\sup_{l\rightarrow+\infty}\left(  \left(  q_{c}+q(l)\right)  /l\right)
<1$). Thus $\mathcal{F}_{\mathcal{L}_{0}}\left(  \mathrm{C}^{\infty}\left(
\Omega\right)  \right)  $ (\textit{resp.} $\mathcal{F}_{\mathcal{L}_{1}%
}\left(  \mathrm{C}^{\infty}\left(  \Omega\right)  \right)  $) are submodulus
of $\mathcal{F}\left(  \mathrm{C}^{\infty}\left(  \Omega\right)  \right)  $.

For $\left(  f_{\varepsilon}\right)  _{\varepsilon}$ and $\left(
g_{\varepsilon}\right)  _{\varepsilon}$) in $\mathcal{F}_{\mathcal{L}_{0}%
}\left(  \mathrm{C}^{\infty}\left(  \Omega\right)  \right)  $, there exists
$C>0$ such that
\[
p_{K,l}\left(  f_{\varepsilon}g_{\varepsilon}\right)  \leq Cp_{K,l}\left(
f_{\varepsilon}\right)  p_{K,l}\left(  g_{\varepsilon}\right)  .
\]
Consequently, $p_{K,l}\left(  f_{\varepsilon}g_{\varepsilon}\right)
=\mathrm{O}\left(  \varepsilon^{-2q(l)}\right)  $ for\ $\varepsilon
\rightarrow0$, with $\lim_{l\rightarrow+\infty}\left(  2q(l)/l\right)  =0$.
Thus $\mathcal{F}_{\mathcal{L}_{0}}\left(  \mathrm{C}^{\infty}\left(
\Omega\right)  \right)  $ is a subalgebra of $\mathcal{F}\left(
\mathrm{C}^{\infty}\left(  \Omega\right)  \right)  $.\medskip\smallskip
\end{proof}

Consequently, we can consider the following subalgebra (\textit{resp}.
submodulus)
\[
\mathcal{G}_{\mathcal{L}_{0}}\left(  \Omega\right)  =\mathcal{F}%
_{\mathcal{L}_{0}}\left(  \mathrm{C}^{\infty}\left(  \Omega\right)  \right)
/\mathcal{N}\left(  \mathrm{C}^{\infty}\left(  \Omega\right)  \right)
\;\;\text{(\textit{resp}. }\mathcal{G}_{\mathcal{L}_{1}}\left(  \Omega\right)
=\mathcal{F}_{\mathcal{L}_{1}}\left(  \mathrm{C}^{\infty}\left(
\Omega\right)  \right)  /\mathcal{N}\left(  \mathrm{C}^{\infty}\left(
\Omega\right)  \,\right)
\]
of $\mathcal{G}\left(  \Omega\right)  $.

\begin{remark}
Some spaces with more restrictive conditions have already been considered (See
e.g. \cite{Ober1}, \cite{Scarpa1}).\ Set
\[
\mathcal{F}^{\infty}\left(  \mathrm{C}^{\infty}\left(  \Omega\right)  \right)
=\left\{  \left(  f_{\varepsilon}\right)  _{\varepsilon}\in\mathcal{F}\left(
\Omega\right)  ^{\left(  0,1\right]  }\,\left|  \forall K\Subset
\Omega,\;\exists q\in\mathbb{N},\;\forall l\in\mathbb{N},\;p_{K,l}\left(
f_{\varepsilon}\right)  =\mathrm{O}\left(  \varepsilon^{-q}\right)
\;\mathrm{for}\;\varepsilon\rightarrow0\right.  \right\}  .
\]
$\mathcal{F}^{\infty}\left(  \mathrm{C}^{\infty}\left(  \Omega\right)
\right)  $ turns to be a subalgebra of $\mathcal{F}_{\mathcal{L}_{0}}\left(
\mathrm{C}^{\infty}\left(  \Omega\right)  \right)  $, $\mathcal{F}%
_{\mathcal{L}_{1}}\left(  \mathrm{C}^{\infty}\left(  \Omega\right)  \right)  $
and
\[
\mathcal{G}^{\infty}\left(  \Omega\right)  =\mathcal{F}^{\infty}\left(
\mathrm{C}^{\infty}\left(  \Omega\right)  \right)  /\mathcal{N}\left(
\mathrm{C}^{\infty}\left(  \Omega\right)  \right)
\]
a subalgebra of $\mathcal{G}_{\mathcal{L}}\left(  \Omega\right)  $,
$\mathcal{G}_{\mathcal{L}_{1}}\left(  \Omega\right)  $ and $\mathcal{G}\left(
\Omega\right)  $. For the local analysis or microlocal analysis of generalized
functions, the $\mathcal{G}^{\infty}$ regularity plays the\ role of the
$\mathrm{C}^{\infty}$' one for distributions \cite{NePiSc} \cite{OberPilSca}.
Our spaces $\mathcal{G}_{\mathcal{L}_{0}}\left(  \Omega\right)  $ and
$\mathcal{G}_{\mathcal{L}_{1},C}\left(  \Omega\right)  $ give new types of
regularity for generalized functions. This will be studied in a forthcoming paper.
\end{remark}

\begin{notation}
We shall note $\mathcal{G}_{C}^{\infty}\left(  \Omega\right)  $ (resp.
$\mathcal{G}_{\mathcal{L}_{0},C}\left(  \Omega\right)  $, $\mathcal{G}%
_{\mathcal{L}_{1},C}\left(  \Omega\right)  $) the subspace of compactly
supported elements of $\mathcal{G}^{\infty}\left(  \Omega\right)  $ (resp.
$\mathcal{G}_{\mathcal{L}_{0}}\left(  \Omega\right)  $, $\mathcal{G}%
_{\mathcal{L}_{1}}\left(  \Omega\right)  $).
\end{notation}

\subsection{Fundamental lemma}

\begin{lemma}
\label{LmnMTh4}Let $d$ be an integer and $\left(  \theta_{\varepsilon}\right)
_{\varepsilon}\in\mathcal{D}\left(  \mathbb{R}^{d}\right)  ^{\left(
0,1\right]  }$ a net of mollifiers satisfying conditions (\ref{LmnGoodM1}) and
(\ref{LmnGoodM2}).\ For any\textbf{\ }$\left(  g_{\varepsilon}\right)
_{\varepsilon}\in\mathcal{F}_{\mathcal{L}_{1}}\left(  \mathrm{C}^{\infty
}\left(  \mathbb{R}^{d}\right)  \right)  $ we have
\begin{equation}
\left(  g_{\varepsilon}\ast\theta_{\varepsilon}-g_{\varepsilon}\right)
_{\varepsilon}\in\mathcal{N}\left(  C^{\infty}\left(  \mathbb{R}^{d}\right)
\right)  . \label{GSTConvolEqu}%
\end{equation}
\end{lemma}

\begin{proof}
We shall prove this lemma in the case $d=1$, the general case only differs by
more complicate algebraic expressions.

Fix $\left(  g_{\varepsilon}\right)  _{\varepsilon}\in F_{\mathcal{L}}\left(
\mathrm{C}^{\infty}\left(  \mathbb{R}^{d}\right)  \right)  $, $K$ a compact of
$\mathbb{R}$ and set $\Delta_{\varepsilon}=g_{\varepsilon}\ast\theta
_{\varepsilon}-g_{\varepsilon}$ for $\varepsilon\in\left(  0,1\right]  $.
Writing $\int\theta_{\varepsilon}\left(  x\right)  \,\mathrm{d}x=1+\mathcal{N}%
_{\varepsilon}$ with $\left(  \mathcal{N}_{\varepsilon}\right)  _{\varepsilon
}\in\mathcal{N}\left(  \mathbb{R}\right)  $ we get
\[
\Delta_{\varepsilon}(y)=\int g_{\varepsilon}(y-x)\theta_{\varepsilon
}(x)\,\mathrm{d}x-g_{\varepsilon}(y)=\int\left(  g_{\varepsilon}%
(y-x)-g_{\varepsilon}(y)\right)  \theta_{\varepsilon}(x)\,\mathrm{d}%
x+\mathcal{N}_{\varepsilon}g_{\varepsilon}(y).
\]
The integration is performed on the compact set $\operatorname*{supp}%
\theta_{\varepsilon}=\left[  -2/\left|  \ln\varepsilon\right|  ,2/\left|
\ln\varepsilon\right|  \right]  $.

Let $m$ be an integer. For each $i\in\mathbb{N},$there exists an integer
$q(i)$ such that
\[
\sup_{\xi\in K^{\prime}}\left|  g_{\varepsilon}^{\left(  i\right)  }\left(
\xi\right)  \right|  =\mathrm{O}\left(  \varepsilon^{-q(i)}\right)
\;\;\mathrm{for}\;\varepsilon\rightarrow0\text{ }%
\]
with $\lim\!\sup_{i\rightarrow+\infty}\left(  q\left(  i\right)  /i\right)  <1
$ and $K^{\prime}$ is a compact such that $\left[  y-1,y+1\right]  \subset
K^{\prime}$ for all $y\in K$.

As $\lim\!\sup_{i\rightarrow+\infty}\left(  q\left(  i\right)  /i\right)  <1$,
we get $\lim_{i\rightarrow+\infty}\left(  i-l(i)\right)  =+\infty$, and there
exists an integer $k$ such that $k-l(k)>m$. Taylor's formula gives
\[
g_{\varepsilon}(y-x)-g_{\varepsilon}(y)=\sum_{i=1}^{k-1}\frac{\left(
-x\right)  ^{i}}{i!}g_{\varepsilon}^{\left(  i\right)  }\left(  y\right)
+\frac{\left(  -x\right)  ^{k-1}}{\left(  k-1\right)  !}\int_{0}%
^{1}g_{\varepsilon}^{(k)}\left(  y-ux\right)  \left(  1-u\right)
^{k-1}\,\mathrm{d}u
\]
and
\begin{multline*}
\Delta_{\varepsilon}(y)=\underset{P_{\varepsilon}\left(  k,y\right)
}{\underbrace{\sum_{i=1}^{k-1}\frac{\left(  -1\right)  ^{i}}{i!}%
g_{\varepsilon}^{\left(  i\right)  }\left(  y\right)  \int x^{i}%
\theta_{\varepsilon}(x)\,\mathrm{d}x}}\\
+\underset{R_{\varepsilon}(k,y)}{\underbrace{\int_{-2/\left|  \ln
\varepsilon\right|  }^{2/\left|  \ln\varepsilon\right|  }\frac{\left(
-x\right)  ^{k-1}}{\left(  k-1\right)  !}\int_{0}^{1}g_{\varepsilon}%
^{(k)}\left(  y-ux\right)  \left(  1-u\right)  ^{k-1}\,\mathrm{d}%
u\,\theta_{\varepsilon}(x)\,\mathrm{d}x}}+\mathcal{N}_{\varepsilon
}g_{\varepsilon}^{(k)}(y).
\end{multline*}
According to lemma \ref{LmnGoodM}, we have $\left(  \int x^{i}\theta
_{\varepsilon}(x)\,\mathrm{d}x\right)  _{\varepsilon}\in\mathcal{N}\left(
\mathbb{R}\right)  $ and consequently
\[
\forall i\in\left\{  0,\ldots,k-1\right\}  ,\;\;\;\int x^{i}\theta
_{\varepsilon}(x)\,\mathrm{d}x=\mathrm{O}\left(  \varepsilon^{m+q(i)}\right)
\;\;\mathrm{for}\;\varepsilon\rightarrow0.
\]
We get
\[
P_{\varepsilon}\left(  k,y\right)  =\mathrm{O}\left(  \varepsilon^{m}\right)
\;\;\mathrm{for}\;\varepsilon\rightarrow0.
\]
Using the definition of $\theta_{\varepsilon}$, we have
\[
R_{\varepsilon}(k,y)=\frac{1}{\varepsilon}\int_{-2/\left|  \ln\varepsilon
\right|  }^{2/\left|  \ln\varepsilon\right|  }\frac{\left(  -x\right)  ^{k-1}%
}{\left(  k-1\right)  !}\left(  \int_{0}^{1}g_{\varepsilon}^{(k)}\left(
y-ux\right)  \left(  1-u\right)  ^{k-1}\,\mathrm{d}u\right)  \,\rho\left(
\frac{x}{\varepsilon}\right)  \chi\left(  x\left|  \ln\varepsilon\right|
\right)  \,\mathrm{d}x.
\]
Setting $v=x/\varepsilon$ we get
\[
R_{\varepsilon}(k,y)=\frac{\varepsilon^{k-1}}{\left(  k-1\right)  !}%
\int_{-2/\left(  \varepsilon\left|  \ln\varepsilon\right|  \right)
}^{2/\left(  \varepsilon\left|  \ln\varepsilon\right|  \right)  }\left(
-v\right)  ^{k-1}\left(  \int_{0}^{1}g_{\varepsilon}^{(k)}\left(
y-\varepsilon uv\right)  \left(  1-u\right)  ^{k-1}\,\mathrm{d}u\right)
\,\rho\left(  v\right)  \chi\left(  \varepsilon\left|  \ln\varepsilon\right|
v\right)  \,\mathrm{d}v.
\]
For $\left(  u,v\right)  \in\left[  0,1\right]  \times\left[  -2/\left(
\varepsilon\left|  \ln\varepsilon\right|  \right)  ,2/\left(  \varepsilon
\left|  \ln\varepsilon\right|  \right)  \right]  $, we have $y-\varepsilon
uv\in\left[  y-1,y+1\right]  $ for $\varepsilon$ small enough. Then, for $y\in
K$, $y-\varepsilon uv$ lies in a compact $K^{\prime}$ for $\left(  u,v\right)
$ in the domain of integration.

It follows
\begin{align*}
\left|  R_{\varepsilon}(k,y)\right|   &  \leq\frac{\varepsilon^{k-1}}{\left(
k-1\right)  !}\sup_{\xi\in K^{\prime}}\left|  g_{\varepsilon}^{(k)}\left(
\xi\right)  \right|  \int_{-2/\left(  \varepsilon\left|  \ln\varepsilon
\right|  \right)  }^{2/\left(  \varepsilon\left|  \ln\varepsilon\right|
\right)  }\left|  v\right|  ^{k-1}\left|  \rho\left(  v\right)  \right|
\mathrm{d}v,\\
&  \leq\frac{\varepsilon^{k-1}}{\left(  k-1\right)  !}\sup_{\xi\in K^{\prime}%
}\left|  g_{\varepsilon}^{(k)}\left(  \xi\right)  \right|  \int_{-\infty
}^{+\infty}\left|  v\right|  ^{k-1}\left|  \rho\left(  v\right)  \right|
\mathrm{d}v,\\
&  \leq C\sup_{\xi\in K^{\prime}}\left|  g_{\varepsilon}^{(k)}\left(
\xi\right)  \right|  \varepsilon^{k-1}\;\;(C>0)
\end{align*}
The constant $C$ depends only on the integer $k$ and $\rho$. By assumption on
$k$, we get
\[
\sup_{y\in K}\left|  R_{\varepsilon}(k,y)\right|  =\mathrm{O}\left(
\varepsilon^{m}\right)  \;\mathrm{for}\;\varepsilon\rightarrow0.
\]
Summering all results, we get $\sup_{y\in K}\Delta_{\varepsilon}%
(y)=\mathrm{O}\left(  \varepsilon^{m}\right)  \;$for$\;\varepsilon\rightarrow0.$

As $\left(  \Delta_{\varepsilon}\right)  _{\varepsilon}\in\mathcal{F}\left(
\mathrm{C}^{\infty}\left(  \mathbb{R}^{d}\right)  \right)  $ and $\sup_{y\in
K}\Delta_{\varepsilon}(y)=\mathrm{O}\left(  \varepsilon^{m}\right)  $
for$\;\varepsilon\rightarrow0$, for all $m>0$ and $K\Subset\mathbb{R}$, we can
conclude that $\left(  \Delta_{\varepsilon}\right)  _{\varepsilon}%
\in\mathcal{N}\left(  \mathrm{C}^{\infty}\left(  \mathbb{R}^{d}\right)
\right)  $ without estimating the derivatives by using theorem 1.2.3. of
\cite{GKOS}.
\end{proof}

\begin{remark}
\label{GSTRemKernIdent}Let us fix a net of mollifiers $\left(  \theta
_{\varepsilon}\right)  _{\varepsilon}$ satisfying conditions (\ref{LmnGoodM1})
and (\ref{LmnGoodM2}) to embed $\mathcal{D}^{\prime}\left(  \mathbb{R}%
^{d}\right)  $ in $\mathcal{G}\left(  \mathbb{R}^{d}\right)  $.\ Relation
(\ref{GSTConvolEqu}) shows that $\left[  \left(  \theta_{\varepsilon}\right)
_{\varepsilon}\right]  $ is plays the role of identity for convolution in
$\mathcal{G}_{\mathcal{L}_{0}}\left(  \mathbb{R}^{d}\right)  $ and
$\mathcal{G}_{\mathcal{L}_{1}}\left(  \mathbb{R}^{d}\right)  $, whereas this
is not true for $\mathcal{G}\left(  \mathbb{R}^{d}\right)  $. This is an
essential feature of these new spaces. (See also example \ref{GSTExamMainTh1} below.)
\end{remark}

\section{Schwartz type theorem\label{GSTSecST}}

\subsection{Extension of linear maps\label{SSBExtension}}

Nets of maps $\left(  L_{\varepsilon}\right)  _{\varepsilon}$ between two
topological algebras having some good growth properties with respect to the
parameter $\varepsilon$ can be canonically extended to the respective
Colombeau spaces based on algebras as it is shown in \cite{ADDS},
\cite{DHPV3}, \cite{GKOS} for examples. We are going to introduce here some
new notions.

We uses the notations of \ref{GSTSSCCompSupp}, specially
\[
\mathcal{D}_{J}\left(  \mathbb{R}^{n}\right)  =\left\{  f\in\mathcal{D}\left(
\mathbb{R}^{n}\right)  \,\left|  \,\operatorname*{supp}f\subset K_{J}\right.
\right\}  ,
\]
where $\left(  K_{J}\right)  _{J\in\mathbb{N}}$ is a sequence of compacts
exhausting $\mathbb{R}^{n}$, and $\mathcal{D}_{J}\left(  \mathbb{R}%
^{n}\right)  $ is endowed with the family of semi norms $p_{J,l}\left(
f\right)  =\sup_{\left|  \alpha\right|  \leq l,\;x\in K_{J}}\left|
\partial^{\alpha}f\left(  x\right)  \right|  .$

\begin{definition}
\label{DefExt1}Let $J$ be an integer and $\left(  L_{\varepsilon}\right)
_{\varepsilon}\in\mathcal{L}\left(  \mathcal{D}_{J}\left(  \mathbb{R}%
^{n}\right)  ,\mathrm{C}^{\infty}\left(  \mathbb{R}^{m}\right)  \right)
^{\left(  0,1\right]  }$ be a net of linear maps.\newline $i.$~We say that
$\left(  L_{\varepsilon}\right)  _{\varepsilon}$ is\emph{\ moderate} if
\begin{multline*}
\forall K\Subset\mathbb{R}^{m},\;\;\forall l\in\mathbb{N},\;\;\exists\left(
C_{\varepsilon}\right)  _{\varepsilon}\in\mathcal{F}\left(  \mathbb{R}%
_{+}\right)  ,\;\;\exists l^{\prime}\in\mathbb{N},\\
\forall f\in\mathcal{D}_{J}\left(  \mathbb{R}^{n}\right)  \;\;\;\;p_{K,l}%
\left(  L_{\varepsilon}\left(  f\right)  \right)  \leq C_{\varepsilon
}p_{J,l^{\prime}}\left(  f\right)  \text{ (for }\varepsilon\text{ small
enough)}.
\end{multline*}
$ii.$~We say that $\left(  L_{\varepsilon}\right)  _{\varepsilon}$ is
\emph{strongly moderate} if
\begin{align*}
\forall K  &  \Subset\mathbb{R}^{m},\;\;\exists\lambda\in\mathbb{N}%
^{\mathbb{N}}\text{ with }\lambda(l)=\mathrm{O}(l)\;\text{for }l\rightarrow
+\infty,\;\;\exists r\in\mathbb{N}^{\mathbb{N}}\text{ with }\underset
{l\rightarrow+\infty}{\lim\sup}\left(  r(l)/l\right)  <1,\\
\forall l  &  \in\mathbb{N},\;\;\exists C\in\mathbb{R}_{+},\;\;\forall
f\in\mathcal{D}_{J}\left(  \mathbb{R}^{n}\right)  ,\;\;\;p_{K,l}\left(
L_{\varepsilon}\left(  f\right)  \right)  \leq C\varepsilon^{-r(l)}%
p_{J,\lambda(l)}\left(  f\right)  \text{ (for }\varepsilon\text{ small
enough)}.
\end{align*}
\end{definition}

In the strong moderation, the growth of $p_{K,l}\left(  L_{\varepsilon}\left(
f\right)  \right)  $ with respect to the index $l$ is controlled by the
sequence $\lambda\left(  \cdot\right)  $ which grows at most like $l$. and by
the sequence $r(l)$.\medskip

As our main result is based on linear maps from $\mathcal{D}\left(
\mathbb{R}^{n}\right)  $ to $\mathrm{C}^{\infty}\left(  \mathbb{R}^{m}\right)
$ we need one further extension:

\begin{definition}
\label{DefExtLimInd}A net of maps $\left(  L_{\varepsilon}\right)
_{\varepsilon}\in$ $\left(  \mathcal{L}\left(  \mathcal{D}\left(
\mathbb{R}^{n}\right)  ,\mathrm{C}^{\infty}\left(  \mathbb{R}^{m}\right)
\right)  ^{\left(  0,1\right]  }\right)  $ is \emph{moderate} (resp.
\emph{strongly moderate}) if for every $J\in\mathbb{N}$, the restriction
$\left(  L_{\varepsilon\left|  \mathcal{D}_{J}\left(  \mathbb{R}^{n}\right)
\right.  }\right)  \in\mathcal{L}\left(  \mathcal{D}_{J}\left(  \mathbb{R}%
^{n}\right)  ,\mathrm{C}^{\infty}\left(  \mathbb{R}^{m}\right)  \right)
^{\left(  0,1\right]  }$ is moderate (resp. strongly moderate) in the sense of
definition \ref{DefExt1}.
\end{definition}

\begin{proposition}
\label{PropExt}Any moderate net $\left(  L_{\varepsilon}\right)  \in$ $\left(
\mathcal{L}\left(  \mathcal{D}\left(  \mathbb{R}^{n}\right)  ,\mathrm{C}%
^{\infty}\left(  \mathbb{R}^{m}\right)  \right)  ^{\left(  0,1\right]
}\right)  $, in the sense of definition \ref{DefExtLimInd}, admits a canonical
extension $L\in\mathcal{L}\left(  \mathcal{G}_{C}\left(  \mathbb{R}%
^{n}\right)  ,\mathcal{G}\left(  \mathbb{R}^{m}\right)  \right)  $ defined by
\begin{equation}
L\left(  \left[  \left(  f_{\varepsilon}\right)  \right]  \right)
=L_{\varepsilon}\left(  f_{\varepsilon}\right)  +\mathcal{N}\left(
\mathrm{C}^{\infty}\left(  \mathbb{R}^{m}\right)  \right)  . \label{GSTextDef}%
\end{equation}
Moreover, if the net $\left(  L_{\varepsilon}\right)  $ is strongly moderate,
$L\left(  \mathcal{G}_{\mathcal{L}_{0},C}\left(  \mathbb{R}^{n}\right)
\right)  $ is included in $\mathcal{G}_{\mathcal{L}_{1}}\left(  \mathbb{R}%
^{m}\right)  $
\end{proposition}

\begin{proof}
Fix $K\Subset\mathbb{R}^{m}$, $l\in\mathbb{N}$ and let $\left(  f_{\varepsilon
}\right)  _{\varepsilon}$ be in $\mathcal{F}_{\mathcal{D}}\left(
\mathbb{R}^{n}\right)  $. There exists $J\in\mathbb{N}$ such that\ $\left(
f_{\varepsilon}\right)  _{\varepsilon}\in\mathcal{F}_{J}\left(  \mathbb{R}%
^{n}\right)  $ and according to the definition of moderate nets, we get
$\left(  C_{\varepsilon}\right)  _{\varepsilon}\in\mathcal{F}\left(
\mathbb{R}_{+}\right)  $ and $l^{\prime}\in\mathbb{N}$ such that
\begin{equation}
p_{K,l}\left(  L_{\varepsilon}\left(  f_{\varepsilon}\right)  \right)  \leq
C_{\varepsilon}p_{J,l^{\prime}}\left(  f_{\varepsilon}\right)  \text{, for
}\varepsilon\text{ small enough}. \label{GSTPrExt1}%
\end{equation}
Inequality (\ref{GSTPrExt1}) leads to $\left(  L_{\varepsilon}\left(
f_{\varepsilon}\right)  \right)  _{\varepsilon}\in\mathcal{F}\left(
\mathrm{C}^{\infty}\left(  \mathbb{R}^{m}\right)  \right)  $. Moreover, if
$\left(  f_{\varepsilon}\right)  _{\varepsilon}$ belongs to $\mathcal{N}%
_{\mathcal{D}}\left(  \mathbb{R}^{n}\right)  $ the same inequality implies
that $\left(  L_{\varepsilon}\left(  f_{\varepsilon}\right)  \right)
_{\varepsilon}\in\mathcal{N}\left(  \mathrm{C}^{\infty}\left(  \mathbb{R}%
^{m}\right)  \right)  $. Those two properties shows that $L$ is well defined
by formula (\ref{GSTextDef}).

Now, suppose that $\left(  L_{\varepsilon}\right)  _{\varepsilon}$ is strongly
moderate and consider $\left(  f_{\varepsilon}\right)  _{\varepsilon}%
\in\mathcal{F}_{\mathcal{L}_{0}}\left(  \mathrm{C}^{\infty}\left(
\mathbb{R}^{n}\right)  \right)  \cap\mathcal{F}_{J}\left(  \mathbb{R}%
^{n}\right)  $. Fix $K\Subset\mathbb{R}^{m}$. There exists a sequence
$\lambda\in\mathbb{N}^{\mathbb{N}}$, with $\lambda(l)=\mathrm{O}(l)\;$for
$l\rightarrow+\infty$, and a sequence $r\in\mathbb{N}^{\mathbb{N}}$ with
$\lim\!\sup_{l\rightarrow+\infty}\left(  r(l)/l\right)  <1\;$such that
\[
\forall l\in\mathbb{N},\;\;\exists C\in\mathbb{R}_{+},\;\;\;p_{K,l}\left(
L_{\varepsilon}\left(  f_{\varepsilon}\right)  \right)  \leq C\varepsilon
^{-r(l)}p_{J,\lambda(l)}\left(  f_{\varepsilon}\right)  \text{ (for
}\varepsilon\text{ small enough)}.
\]
As $\left(  f_{\varepsilon}\right)  _{\varepsilon}$ is in $\mathcal{F}%
_{\mathcal{L}_{0}}\left(  \mathrm{C}^{\infty}\left(  \mathbb{R}^{n}\right)
\right)  $, there exists a sequence $q\in\mathbb{N}^{\mathbb{N}},\;$with
$\lim_{\lambda\rightarrow+\infty}\left(  q(\lambda)/\lambda\right)  =0$ such
that
\[
\forall\lambda\in\mathbb{N},\;\;\;\;\;p_{J,\lambda}\left(  f_{\varepsilon
}\right)  =\mathrm{O}\left(  \varepsilon^{-q(\lambda)}\right)  \;\mathrm{for}%
\;\varepsilon\rightarrow0.
\]
We get that
\[
\forall l\in\mathbb{N},\;\;\;\;\;p_{K,l}\left(  L_{\varepsilon}\left(
f_{\varepsilon}\right)  \right)  =\mathrm{O}\left(  \varepsilon^{-q_{1}\left(
l\right)  }\right)  \;\;\mathrm{for}\;\varepsilon\rightarrow0,\;\text{\ with
}q_{1}\left(  l\right)  =r(l)+q\left(  \lambda(l)\right)  .
\]
If $\lambda(l)$ is bounded, we get immediately that $q_{1}\left(  l\right)
/l=\mathrm{o}\left(  1\right)  $ for $l\rightarrow+\infty$. If $\lambda(l)$ is
not bounded, we write for $l$ such that $\lambda(l)\neq0$.
\[
\frac{q_{1}\left(  l\right)  }{l}=\frac{r(l)}{l}+\frac{q\left(  \lambda
(l)\right)  }{\lambda(l)}\frac{\lambda(l)}{l}%
\]
Since $\lambda(l)/l$ is bounded and $q\left(  l\right)  /l=\mathrm{o}\left(
1\right)  $, we get that $\frac{q\left(  \lambda(l)\right)  }{\lambda(l)}%
\frac{\lambda(l)}{l}=\mathrm{o}\left(  1\right)  $ for $l\rightarrow+\infty$.
This gives that $\lim\!\sup_{l\rightarrow+\infty}\left(  q_{1}(l)/l\right)
<1$ and $\left(  L_{\varepsilon}\left(  f_{\varepsilon}\right)  \right)
_{\varepsilon}\in\mathcal{F}_{\mathcal{L}_{1}}\left(  \mathrm{C}^{\infty
}\left(  \Omega\right)  \right)  $ and shows last assertion.
\end{proof}

\subsection{Main theorem}

\begin{theorem}
\label{ThmGST}Let $\left(  L_{\varepsilon}\right)  _{\varepsilon}%
\in\mathcal{L}\left(  \mathcal{D}(\mathbb{R}^{n}),\mathrm{C}^{\infty
}(\mathbb{R}^{m})\right)  ^{\left(  0,1\right]  }$ be a net of strongly
moderate continuous linear maps and $L\in\mathcal{L}\left(  \mathcal{G}%
_{C}(\mathbb{R}^{n}),\mathcal{G}(\mathbb{R}^{m})\right)  ^{\left(  0,1\right]
}$ its canonical extension. There exists $H_{L}\in\mathcal{G}\left(
\mathbb{R}^{m}\times\mathbb{R}^{n}\right)  $ such that
\[
\forall f\in\mathcal{G}_{\mathcal{L}_{_{0}},C}\left(  \mathbb{R}^{n}\right)
),\;\;\;\;L\left(  f\right)  (x)=\int H_{L}\left(  x,y\right)
f(y)\,\mathrm{d}y.
\]
\end{theorem}

In other words, $L$ restricted to $\mathcal{G}_{\mathcal{L}_{0},C}\left(
\Omega\right)  )$ can be represented by a kernel $H_{L}$. The fact that the
equality is only valid in $\mathcal{G}_{\mathcal{L}_{0},C}\left(
\mathbb{R}^{n}\right)  )$ is not surprising. The structure of the theorem is
similar as Schwartz'one: $f$ belongs to a ``smaller'' type of space as $H_{L}$
and $L\left(  f\right)  $, which both belongs to the same kind of space.

\begin{example}
\label{GSTExamMainTh1}Remark \ref{GSTRemKernIdent} and relation
(\ref{GSTConvolEqu}) shows also that the identity map of $\mathcal{G}%
_{\mathcal{L}_{0},C}\left(  \mathbb{R}^{n}\right)  )$ admits as kernel
\begin{equation}
\Phi=cl\left(  \left(  \left(  x,y\right)  \mapsto\varphi_{\varepsilon}\left(
x-y\right)  \right)  _{\varepsilon}\right)  \label{GSTExample1}%
\end{equation}
where $\left(  \varphi_{\varepsilon}\right)  _{\varepsilon\in\left(
0,1\right]  }$ is any net of mollifiers satisfying conditions (\ref{LmnGoodM1}%
) and (\ref{LmnGoodM2}) of lemma \ref{LmnGoodM}.
\end{example}

This example shows also that, in general, we don't have uniqueness in theorem
\ref{ThmGST}, but a so called \emph{weak uniqueness}. In our example, any net
$\left(  \varphi_{\varepsilon}\right)  _{\varepsilon}$ of mollifiers satisfies
$\varphi_{\varepsilon}\rightarrow\delta$ in $\mathcal{D}^{\prime}$ for
$\varepsilon\rightarrow0$: Thus, kernels of the form (\ref{GSTExample1}) are
associated in $\mathcal{G}\left(  \mathbb{R}^{m}\times\mathbb{R}^{n}\right)  $
or weakly equal i.e. the difference of their representative tends to $0$ in
$\mathcal{D}^{\prime}$ for $\varepsilon\rightarrow0$. (See \cite{DHPV3},
\cite{GKOS}, \cite{JAM2}, \cite{NePiSc} for further analysis of different
associations in Colombeau type spaces.)

\subsection{Link with the classical Schwartz theorem: Equality in generalized
distribution sense\label{GSTSecRship}}

Let $\Lambda\in\mathcal{L}\left(  \mathcal{D}\left(  \mathbb{R}^{n}\right)
,\mathcal{D}^{\prime}\left(  \mathbb{R}^{m}\right)  \right)  $ be continuous
for the strong topology and consider the family of linear mappings $\left(
L_{\varepsilon}\right)  _{\varepsilon\in}$ defined by
\[
L_{\varepsilon}:\mathcal{D}\left(  \mathbb{R}^{n}\right)  \rightarrow
\mathrm{C}^{\infty}\left(  \mathbb{R}^{m}\right)  \;\;\;\;f\mapsto
\Lambda\left(  f\right)  \ast\varphi_{\sqrt{\varepsilon}}\text{,}%
\]
where $\left(  \varphi_{\varepsilon}\right)  _{\varepsilon}$ is a family of
mollifiers satisfying conditions (\ref{LmnGoodM1}) and (\ref{LmnGoodM2}) of
lemma \ref{LmnGoodM}. We have:

\begin{proposition}
\label{GSTCoroTHS}\ \ \newline $i$.~For all $\varepsilon\in\left(  0,1\right]
$, $L_{\varepsilon}$ is continuous for the usual topologies of $\mathcal{D}%
\left(  \mathbb{R}^{n}\right)  $ and $\mathrm{C}^{\infty}\left(
\mathbb{R}^{m}\right)  $.\ \newline $ii$.~The net $\left(  L_{\varepsilon
}\right)  _{\varepsilon}$ is strongly moderate.
\end{proposition}

Consequently, theorem \ref{ThmGST} shows that the canonical extension $L$ of
the net $\left(  L_{\varepsilon}\right)  _{\varepsilon}$ admits a kernel
$H_{L}$.

\begin{proposition}
\label{GSTEequaGD}For all $f\in\mathcal{D}\left(  \mathbb{R}^{n}\right)  $,
$\Lambda\left(  f\right)  $ is equal to $\widetilde{H}_{L}\left(  f\right)  $
in the \emph{generalized distribution sense}, that is
\[
\forall\Phi\in\mathcal{D}\left(  \mathbb{R}^{m}\right)  ,\;\;\;\;\left\langle
\Lambda\left(  f\right)  ,\Phi\right\rangle =\left\langle \widetilde{H}%
_{L}\left(  f\right)  ,\Phi\right\rangle \text{ in }\overline{\mathbb{C}%
}\text{.}%
\]
\end{proposition}

This generalized distribution equality, introduced in \cite{NePiSc}, means in
other words that, for all $k\in\mathbb{N}$,
\begin{equation}
\forall\Phi\in\mathcal{D}\left(  \mathbb{R}^{m}\right)  ,\;\;\left\langle
\Lambda\left(  f\right)  ,\Phi\right\rangle -\int\left(  \int H_{L,\varepsilon
}\left(  x,y\right)  f\left(  y\right)  \,\mathrm{d}y\right)  \Phi\left(
x\right)  \,\mathrm{d}x=\mathrm{O}\left(  \varepsilon^{k}\right)  \text{, for
}\varepsilon\rightarrow0, \label{GSTDefgdequa}%
\end{equation}
where $\left(  H_{L,\varepsilon}\right)  _{\varepsilon}$ is any representative
of $H_{L}$.\smallskip\medskip

In particular, this result implies that $\Lambda\left(  f\right)  $ and
$\widetilde{H}_{L}\left(  f\right)  $ are associated or weakly equal, \emph{id
est}%
\[
\left\{  x\mapsto\int H_{L,\varepsilon}\left(  x,y\right)  f\left(  y\right)
\,\mathrm{d}y\right\}  \longrightarrow\Lambda\left(  f\right)  \text{ in
}\mathcal{D}^{\prime}\text{ for }\varepsilon\rightarrow0\text{.}%
\]

\section{Proofs of theorem \ref{ThmGST} and propositions \ref{GSTCoroTHS} and
\ref{GSTEequaGD}}

\subsection{Proof of theorem \ref{ThmGST}}

Let us fix $\left(  \varphi_{\varepsilon}\right)  _{\varepsilon}\in\left(
\mathcal{D}\left(  \mathbb{R}^{m}\right)  \right)  ^{\left(  0,1\right]  }$
(\emph{resp}. $\left(  \psi_{\varepsilon}\right)  _{\varepsilon}\in\left(
\mathcal{D}\left(  \mathbb{R}^{m}\right)  \right)  ^{\left(  0,1\right]  }$) a
net of mollifiers satisfying conditions \ref{LmnGoodM1} and \ref{LmnGoodM2} of
lemma \ref{LmnGoodM}. For all $y\in\mathbb{R}^{n}$ we define
\[
\psi_{\varepsilon,.}:\mathbb{R}^{n}\rightarrow\mathcal{D}\left(
\mathbb{R}^{n}\right)  \;\;\;\;y\mapsto\psi_{\varepsilon,y}=\left\{
v\mapsto\psi_{\varepsilon}\left(  y-v\right)  \right\}  .
\]

For all $y\in\mathbb{R}^{n}$ and $\varepsilon\in\left(  0,1\right]  $, we set
$\Psi_{\varepsilon,y}=L_{\varepsilon}\left(  \psi_{\varepsilon,y}\right)  $.

\begin{lemma}
\label{LmnMTh1}The map
\[
\Psi_{\varepsilon}:\mathbb{R}^{n}\rightarrow\mathrm{C}^{\infty}(\mathbb{R}%
^{m})\;\;\;\;y\mapsto\Psi_{\varepsilon,y}=L_{\varepsilon}\left(
\psi_{\varepsilon,y}\right)
\]
is of class $\mathrm{C}^{\infty}$ for all $\varepsilon\in\left(  0,1\right]  $.
\end{lemma}

\begin{proof}
The map $\left(  y,v\right)  \mapsto\psi_{\varepsilon}\left(  y-v\right)  $
from $\mathbb{R}^{2n}$ to $\mathbb{R}$ is clearly of class $\mathrm{C}%
^{\infty}$. It follows that the map $\psi_{\varepsilon,\cdot}:y\mapsto
\psi_{\varepsilon,y}$, considered as a map from $\mathbb{R}^{n}$ to
$\mathrm{C}^{\infty}(\mathbb{R}^{n})$, is $\mathrm{C}^{\infty}$ (see for
example theorem 2.2.2 of \cite{GKOS}). As each $\psi_{\varepsilon,y}$ is
compactly supported we can show that $\psi_{\varepsilon,\cdot}$ belongs in
fact to $\mathrm{C}^{\infty}\left(  \mathbb{R}^{n},\mathcal{D}\left(
\mathbb{R}^{n}\right)  \right)  $ by using local arguments. Since
$L_{\varepsilon}$ is linear and continuous it follows that $\Psi_{\varepsilon
}$ is $\mathrm{C}^{\infty}$.\bigskip
\end{proof}

Let us define, for all $\varepsilon\in\left(  0,1\right]  $ and $\left(
x,y\right)  \in\mathbb{R}^{m}\times\mathbb{R}^{n}$:
\[
H_{\varepsilon}\left(  x,y\right)  =\left(  \Psi_{\varepsilon,y}\ast
\varphi_{\varepsilon}\right)  \left(  x\right)  =\int L_{\varepsilon}\left(
\psi_{\varepsilon,y}\right)  \left(  x-\lambda\right)  \varphi_{\varepsilon
}\left(  \lambda\right)  \,\mathrm{d}\lambda.
\]
Note that, for all $\varepsilon\in\left(  0,1\right]  $, this integral is
performed on the compact set $\operatorname*{supp}\varphi_{\varepsilon}$.

\begin{lemma}
\label{LmnMTh2}For all $\varepsilon\in\left(  0,1\right]  $, $H_{\varepsilon}$
is of class $\mathrm{C}^{\infty}$ and $\left(  H_{\varepsilon}\right)
_{\varepsilon}\in\mathcal{F}\left(  \mathbb{R}^{m}\times\mathbb{R}^{n}\right)
$.
\end{lemma}

\begin{proof}
First, the map $g\mapsto g\ast\varphi_{\varepsilon}$ from $\mathrm{C}^{\infty
}(\mathbb{R}^{m})$ into itself is linear continuous and therefore
$\mathrm{C}^{\infty}$. Using lemma \ref{LmnMTh1}, we get that the map
$y\mapsto\left(  \Psi_{\varepsilon,y}\ast\varphi_{\varepsilon}\right)
=H_{\varepsilon}\left(  \cdot,y\right)  $ from $\mathbb{R}^{n}$ to
$\mathrm{C}^{\infty}(\mathbb{R}^{m})$ is $\mathrm{C}^{\infty}$. Using again
theorem 2.2.2 of \cite{GKOS}, we get that $H_{\varepsilon}$ belongs to
$\mathrm{C}^{\infty}(\mathbb{R}^{2n})$.

Consider $K$ and $K^{\prime}$ two compact subsets of $\mathbb{R}^{n}$. Let us
recall that the support of $\psi_{\varepsilon}$ is compact and decreasing to
$\left\{  0\right\}  $ when $\varepsilon$ tends to $0$. Then, there exists a
compact set $K_{\psi}\subset\mathbb{R}^{m}$ such that, for all $\varepsilon
\in\left(  0,1\right]  $, $\operatorname*{supp}\psi_{\varepsilon}\subset
K_{\psi}$ and $\operatorname*{supp}\psi_{\varepsilon,y}\subset y-K\psi$.
Moreover, we can find a compact $K_{J}$ (notation are those of
\ref{SSBExtension}) such that
\[
\forall\varepsilon\in\left(  0,1\right]  ,\;\ \ \;\forall y\in K^{\prime
},\;\;\;\psi_{\varepsilon,y}\in\mathcal{D}_{J}\left(  \mathbb{R}^{n}\right)
.
\]
and $p_{J,l}(\psi_{\varepsilon,y})=p_{K_{\psi},l}(\psi_{\varepsilon})$, for
all $\varepsilon\in\left(  0,1\right]  $.

Let now consider $\left(  \alpha,\beta\right)  \in\left(  \mathbb{N}%
^{n}\right)  ^{2}$ and $\partial^{\alpha}$ (\textit{resp}. $\partial^{\beta}$)
the $\alpha$-partial derivative (\textit{resp}. $\beta$-partial derivative)
with respect to the variable $x$ (\textit{resp}. $y$). Noticing that there
exists a compact set $K_{\varphi}\subset\mathbb{R}^{m}$ such that, for all
$\varepsilon\in\left(  0,1\right]  $, $\operatorname*{supp}\varphi
_{\varepsilon,y}\subset K_{\varphi}$ we get the existence of a constant $C$
such that, for all $\varepsilon\in\left(  0,1\right]  $,
\begin{align*}
\left.  \forall\left(  x,y\right)  \in K\times K^{\prime},\;\;\left|
H_{\varepsilon}\left(  x,y\right)  \right|  \right.   &  \leq C\sup_{\xi\in
K-K\varphi}\left|  \partial^{\beta}L_{\varepsilon}\left(  \psi_{\varepsilon
,y}\right)  \left(  \xi\right)  \right|  \sup_{\xi\in K\varphi}\left|
\partial^{\alpha}\varphi_{\varepsilon}\left(  \xi\right)  \right|  ,\\
&  \leq Cp_{K-K\varphi,\left|  \beta\right|  }\left(  L_{\varepsilon}\left(
\psi_{\varepsilon,y}\right)  \right)  p_{K\varphi,\left|  \alpha\right|
}\left(  \varphi_{\varepsilon}\right)  .
\end{align*}
The moderateness of $\left(  L_{\varepsilon}\right)  _{\varepsilon}$ implies
the existence of $l\in\mathbb{N}$ and $\left(  C_{\varepsilon}^{\prime
}\right)  _{\varepsilon}\in\mathcal{F}\left(  \mathbb{R}_{+}\right)  $ such
that, for all $\varepsilon\in\left(  0,1\right]  $,
\[
\forall\left(  x,y\right)  \in K\times K^{\prime},\;\;\left|  H_{\varepsilon
}\left(  x,y\right)  \right|  \leq C_{\varepsilon}^{\prime}p_{J,l}\left(
\psi_{\varepsilon,y}\right)  p_{K\varphi,\left|  \alpha\right|  }\left(
\varphi_{\varepsilon}\right)  \leq C_{\varepsilon}^{\prime}p_{K_{\psi},l}%
(\psi_{\varepsilon})p_{K\varphi,\left|  \alpha\right|  }\left(  \varphi
_{\varepsilon}\right)  .
\]
The last inequality shows that $\left(  p_{K\times K^{\prime}\left|
\alpha\right|  +\left|  \beta\right|  }(H_{\varepsilon})\right)
_{\varepsilon}$ belongs to $\mathcal{F}\left(  \mathbb{R}_{+}\right)  $, this
ending the proof.\bigskip
\end{proof}

For all $\left(  f_{\varepsilon}\right)  _{\varepsilon}$ in $\mathcal{F}%
\left(  \mathcal{D}\left(  \mathbb{R}^{n}\right)  \right)  $ (this set is
defined by relation (\ref{ModCompSup})) we can consider
\[
\widetilde{H}_{\varepsilon}\left(  f_{\varepsilon}\right)  \left(  x\right)
=\int H_{\varepsilon}\left(  x,y\right)  f_{\varepsilon}(y)\,\mathrm{d}%
y=\int\left(  \int L_{\varepsilon}\left(  \psi_{\varepsilon,y}\right)  \left(
x-\lambda\right)  \varphi_{\varepsilon}\left(  \lambda\right)  \,\mathrm{d}%
\lambda\right)  f_{\varepsilon}(y)\,\mathrm{d}y.
\]
since for all $\varepsilon\in\left(  0,1\right]  $, $f_{\varepsilon}$ is
compactly supported.

\begin{lemma}
\label{LmnMTh3}For all $\left(  f_{\varepsilon}\right)  _{\varepsilon}$ in
$\mathcal{F}\left(  \mathcal{D}\left(  \mathbb{R}^{n}\right)  \right)  $, we
have
\[
\widetilde{H}_{\varepsilon}\left(  f_{\varepsilon}\right)  \left(  x\right)
=\left(  L_{\varepsilon}\left(  \psi_{\varepsilon}\ast f_{\varepsilon}\right)
\ast\varphi_{\varepsilon}\right)  (x).
\]
\end{lemma}

\begin{proof}
Let $\left(  f_{\varepsilon}\right)  _{\varepsilon}$ be in $\mathcal{F}\left(
\mathcal{D}\left(  \mathbb{R}^{n}\right)  \right)  $. For all $\varepsilon
\in\left(  0,1\right]  $ and $x\in\mathbb{R}^{m}$, we have
\begin{align*}
\widetilde{H}_{\varepsilon}\left(  f_{\varepsilon}\right)  \left(  x\right)
&  =\int_{\operatorname*{supp}f}\left(  \int_{\operatorname*{supp}%
\varphi_{\varepsilon}}L_{\varepsilon}\left(  \psi_{\varepsilon,y}\right)
\left(  x-\lambda\right)  \varphi_{\varepsilon}\left(  \lambda\right)
\,\mathrm{d}\lambda\right)  f_{\varepsilon}(y)\,\mathrm{d}y,\\
&  =\int_{\operatorname*{supp}\varphi_{\varepsilon}}\int_{\operatorname*{supp}%
f}L_{\varepsilon}\left(  \psi_{\varepsilon,y}\right)  \left(  x-\lambda
\right)  \varphi_{\varepsilon}\left(  \lambda\right)  f_{\varepsilon
}(y)\,\mathrm{d}\lambda\mathrm{d}y,\\
&  =\int\left(  \int L_{\varepsilon}\left(  \psi_{\varepsilon,y}\right)
\left(  x-\lambda\right)  f_{\varepsilon}(y)\,\mathrm{d}y\right)
\varphi_{\varepsilon}\left(  \lambda\right)  \,\mathrm{d}\lambda,
\end{align*}
the two last equalities being true by Fubini's theorem, each integral being
calculated on a compact set.

For all $\varepsilon\in\left(  0,1\right]  $ and $\xi\in\mathbb{R}^{m}$, we
have the following equality:
\begin{align*}
\int L_{\varepsilon}\left(  \psi_{\varepsilon,y}\right)  \left(  \xi\right)
f_{\varepsilon}(y)\,\mathrm{d}y  &  =L_{\varepsilon}\left(  v\mapsto\int
\psi_{\varepsilon,y}\left(  v\right)  f_{\varepsilon}(y)\mathrm{d}y\right)
\left(  \xi\right)  ,\\
&  =L_{\varepsilon}\left(  v\mapsto\int\psi_{\varepsilon}\left(  y-v\right)
f_{\varepsilon}(y)\mathrm{d}y\right)  \left(  \xi\right)  .
\end{align*}
Indeed, the integrals under consideration in the above equalities are
integrals of continuous functions on compact sets and can be considered as
limits of Riemann sums in the spirit of \cite{HorPDOT1} (Lemma 4.1.3,
p.\ 89):
\begin{gather*}
\forall\xi\in\mathbb{R}^{m},\;\;\int L_{\varepsilon}\left(  \psi
_{\varepsilon,y}\right)  \left(  \xi\right)  f_{\varepsilon}(y)\,\mathrm{d}%
y=\lim_{h\rightarrow0}\sum_{k\in\mathbb{Z}}h^{n}L_{\varepsilon}\left(
\psi_{\varepsilon}\left(  kh-v\right)  \right)  \left(  \xi\right)
f_{\varepsilon}(kh),\\
\forall v\in\mathbb{R}^{n},\;\;\int\psi_{\varepsilon}\left(  y-v\right)
f_{\varepsilon}(y)\mathrm{d}y=\lim_{h\rightarrow0}\sum_{k\in\mathbb{Z}}%
h^{n}\psi_{\varepsilon}\left(  kh-v\right)  f_{\varepsilon}kh).
\end{gather*}
As the mapping $L_{\varepsilon}$ is linear, we have
\[
L_{\varepsilon}\left(  \sum_{k\in\mathbb{Z}}\psi_{\varepsilon}\left(
kh-v\right)  f_{\varepsilon}(kh)\right)  =\sum_{k\in\mathbb{Z}}f_{\varepsilon
}(kh)L_{\varepsilon}\left(  \psi_{\varepsilon}\left(  kh-v\right)  \right)  ,
\]
as each $f_{\varepsilon}(kh)$ is a scalar: The function $\psi_{\varepsilon,y}$
is on the $v$ variable, belonging to $\mathbb{R}^{n}$. By continuity of
$L_{\varepsilon}$, we get
\begin{align*}
L_{\varepsilon}\left(  \int\psi_{\varepsilon}\left(  y-v\right)
f_{\varepsilon}(y)\mathrm{d}y\right)  \left(  \xi\right)   &  =L_{\varepsilon
}\left(  \lim_{h\rightarrow0}\sum_{k\in\mathbb{Z}}h^{n}\psi_{\varepsilon
}\left(  kh-v\right)  f_{\varepsilon}(kh)\right)  \left(  \xi\right)  ,\\
&  =\lim_{h\rightarrow0}\left(  \sum_{k\in\mathbb{Z}}f_{\varepsilon
}(kh)L_{\varepsilon}\left(  \psi_{\varepsilon}\left(  kh-v\right)  \right)
\left(  \xi\right)  \right)  ,\\
&  =\int L_{\varepsilon}\left(  \psi_{\varepsilon,y}\right)  \left(
\xi\right)  f_{\varepsilon}(y)\,\mathrm{d}y.
\end{align*}
Finally, we get for all $\varepsilon\in\left(  0,1\right]  $ and $\xi
\in\mathbb{R}^{m}$,
\begin{align*}
\int L_{\varepsilon}\left(  \psi_{\varepsilon,y}\right)  \left(  \xi\right)
f_{\varepsilon}(y)\,\mathrm{d}y  &  =L_{\varepsilon}\left(  \int
\psi_{\varepsilon}\left(  y-v\right)  f_{\varepsilon}(y)\mathrm{d}y\right)
\left(  \xi\right)  ,\\
&  =L_{\varepsilon}\left(  \psi_{\varepsilon}\ast f_{\varepsilon}\right)
\left(  \xi\right)  ,
\end{align*}
and
\begin{equation}
\widetilde{H}_{\varepsilon}\left(  f_{\varepsilon}\right)  \left(  x\right)
=\int L_{\varepsilon}\left(  \psi_{\varepsilon}\ast f_{\varepsilon}\right)
\left(  x-\lambda\right)  \varphi_{\varepsilon}\left(  \lambda\right)
\,\mathrm{d}\lambda=\left(  L_{\varepsilon}\left(  \psi_{\varepsilon}\ast
f_{\varepsilon}\right)  \ast\varphi_{\varepsilon}\right)  (x).
\label{ThmPValHF}%
\end{equation}
\medskip\smallskip
\end{proof}

We are now complete the proof of theorem \ref{ThmGST}.\ Set
\[
H_{L}=\left(  H_{\varepsilon}\right)  _{\varepsilon}+\mathcal{N}\left(
C^{\infty}\left(  \mathbb{R}^{m+n}\right)  \right)  =\left(  \left(
x,y\right)  \mapsto\left(  \,\Psi_{\varepsilon,y}\ast\varphi_{\varepsilon
}\right)  \left(  x\right)  \,\right)  _{\varepsilon}+\mathcal{N}\left(
C^{\infty}\left(  \mathbb{R}^{m+n}\right)  \right)  .
\]
For all $\left(  f_{\varepsilon}\right)  _{\varepsilon}$ in $\mathcal{F}%
_{\mathcal{L}_{0}}\left(  \mathcal{D}\left(  \mathbb{R}^{n}\right)  \right)  $
we have
\[
\widetilde{H}_{L}\left(  \left[  \left(  f_{\varepsilon}\right)
_{\varepsilon}\right]  \right)  =\left[  \left(  \widetilde{H}_{\varepsilon
}\left(  f_{\varepsilon}\right)  \right)  _{\varepsilon}\right]
\]
by definition of the integral in $\mathcal{G}\left(  \mathbb{R}^{n}\right)
$.\ We have to compare $\left(  \widetilde{H}_{\varepsilon}\left(
f_{\varepsilon}\right)  \right)  _{\varepsilon}$ and $\left(  L_{\varepsilon
}\left(  f_{\varepsilon}\right)  \right)  _{\varepsilon}$. According to lemma
\ref{LmnMTh3}, we have for all $\varepsilon\in\left(  0,1\right]  $,
\begin{align*}
\widetilde{H}_{\varepsilon}\left(  f_{\varepsilon}\right)  -L_{\varepsilon
}\left(  f_{\varepsilon}\right)   &  =\left(  L_{\varepsilon}\left(
\psi_{\varepsilon}\ast f_{\varepsilon}\right)  \ast\varphi_{\varepsilon
}\right)  -L_{\varepsilon}\left(  f_{\varepsilon}\right)  ,\\
&  =L_{\varepsilon}\left(  \psi_{\varepsilon}\ast f_{\varepsilon}\right)
\ast\varphi_{\varepsilon}-L_{\varepsilon}\left(  f_{\varepsilon}\right)
\ast\varphi_{\varepsilon}+L_{\varepsilon}\left(  f_{\varepsilon}\right)
\ast\varphi_{\varepsilon}-L_{\varepsilon}\left(  f_{\varepsilon}\right)  ,\\
&  =L_{\varepsilon}\left(  \psi_{\varepsilon}\ast f_{\varepsilon
}-f_{\varepsilon}\right)  \ast\varphi_{\varepsilon}+L_{\varepsilon}\left(
f_{\varepsilon}\right)  \ast\varphi_{\varepsilon}-L_{\varepsilon}\left(
f_{\varepsilon}\right)  .
\end{align*}
Remarking that $\left(  f_{\varepsilon}\right)  _{\varepsilon}\in
\mathcal{F}_{\mathcal{L}_{0}}\left(  \mathrm{C}^{\infty}\left(  \Omega\right)
\right)  $ and $\left(  L_{\varepsilon}\left(  f_{\varepsilon}\right)
\right)  _{\varepsilon}\in\mathcal{F}_{\mathcal{L}_{1}}\left(  \mathrm{C}%
^{\infty}\left(  \Omega\right)  \right)  $ we get $\left(  \,L_{\varepsilon
}\left(  f_{\varepsilon}\right)  \ast\varphi_{\varepsilon}-L_{\varepsilon
}\left(  f_{\varepsilon}\right)  \,\right)  _{\varepsilon}\in\mathcal{N}%
\left(  \mathrm{C}^{\infty}\left(  \mathbb{R}^{m}\right)  \right)  $ and
$\left(  \,\psi_{\varepsilon}\ast f_{\varepsilon}-f_{\varepsilon}\,\right)
_{\varepsilon}\in\mathcal{N}\left(  \mathrm{C}^{\infty}\left(  \mathbb{R}%
^{m}\right)  \right)  $ by lemma \ref{LmnMTh4}. This last property gives
\[
\left(  \,L_{\varepsilon}\left(  \psi_{\varepsilon}\ast f_{\varepsilon
}-f_{\varepsilon}\right)  \,\right)  _{\varepsilon}\in\mathcal{N}\left(
\mathrm{C}^{\infty}\left(  \mathbb{R}^{m}\right)  \right)  \;\;\mathrm{and}%
\;\;\left(  \,L_{\varepsilon}\left(  \psi_{\varepsilon}\ast f_{\varepsilon
}-f_{\varepsilon}\right)  \ast\varphi_{\varepsilon}\,\right)  \,\in
\mathcal{N}\left(  \mathrm{C}^{\infty}\left(  \mathbb{R}^{m}\right)  \right)
,
\]
since $\left(  \eta_{\varepsilon}\ast\varphi_{\varepsilon}\right)
_{\varepsilon}\in\mathcal{N}\left(  \mathrm{C}^{\infty}\left(  \mathbb{R}%
^{m}\right)  \right)  $ for all $\left(  \eta_{\varepsilon}\right)
_{\varepsilon}\in\mathcal{N}\left(  \mathrm{C}^{\infty}\left(  \mathbb{R}%
^{m}\right)  \right)  $. Finally
\[
\left[  \left(  \widetilde{H}_{\varepsilon}\left(  f_{\varepsilon}\right)
\right)  _{\varepsilon}\right]  =\left[  \left(  L_{\varepsilon}\left(
f_{\varepsilon}\right)  \right)  _{\varepsilon}\right]  =L\left(  \left[
\left(  f_{\varepsilon}\right)  _{\varepsilon}\right]  \right)  ,
\]
this last equality by definition of the extension of a linear map.

\subsection{Proof of proposition \ref{GSTCoroTHS}}

\textit{Assertion }$i$.~We have only to prove continuity on $0$. Let us fix
$\varepsilon\in\left(  0,1\right]  $. Take $\left(  f_{k}\right)  _{k}%
\in\mathcal{D}\left(  \mathbb{R}^{n}\right)  ^{\mathbb{N}}$ a sequence
converging to $0$ in $\mathcal{D}\left(  \mathbb{R}^{n}\right)  $. Since
$\Lambda$ is continuous, the sequence $\left(  T_{k}\right)  _{k}=\left(
\Lambda\left(  f_{k}\right)  \right)  _{k}$ tends to $0$ in $\mathcal{D}%
^{\prime}\left(  \mathbb{R}^{m}\right)  $ for the strong topology. Let us
recall that \cite{Schwartz1}:

\begin{lemma}
\label{LmnSch}A sequence $\left(  T_{k}\right)  _{k}$ tends to $0$ in
$\mathcal{D}^{\prime}\left(  \mathbb{R}^{m}\right)  $ for the strong topology
if, and only if for all $\theta\in\mathcal{D}\left(  \mathbb{R}^{m}\right)  $
the sequence $\left(  T_{k}\ast\theta\right)  _{k}$ tends to $0$, uniformly on
every compact set.
\end{lemma}

For all $\alpha$ in $\mathbb{N}^{m}$, we take $\theta_{\alpha}=\partial
^{\alpha}\varphi_{\sqrt{\varepsilon}}$. Applying lemma \ref{LmnSch}, the
sequences
\[
\left(  T_{k}\ast\partial^{\alpha}\varphi_{\sqrt{\varepsilon}}\right)
_{k}=\left(  \partial^{\alpha}\left(  T_{k}\ast\varphi_{\sqrt{\varepsilon}%
}\right)  \right)  _{k}%
\]
tends to $0$ uniformly on each compact of $\mathbb{R}^{m}$. Then
$L_{\varepsilon}$ is continuous.\medskip

\noindent\textit{Assertion }$ii$.~According to definition \ref{DefExtLimInd},
we have to show that, for all $J\in\mathbb{N}$, the net $\left(
L_{\varepsilon\left|  \mathcal{D}_{J}\right.  }\right)  _{\varepsilon}%
\in\left(  \mathcal{L}\left(  \mathcal{D}_{J}\left(  \mathbb{R}^{n}\right)
,\mathcal{D}^{\prime}\left(  \mathbb{R}^{m}\right)  \right)  \right)
^{\left(  0,1\right]  }$ is strongly moderate. We have
\begin{align*}
\left.  \forall f\in\mathcal{D}_{J}\left(  \mathbb{R}^{n}\right)  ,\;\forall
x\in\mathbb{R}^{m},\;\forall\alpha\in\mathbb{N}^{m},\;\;\partial^{\alpha
}\left(  L_{\varepsilon\left.  \mathcal{D}_{J}\right|  }(f)\right)
(x)\right.   &  =\left(  \Lambda\left(  f\right)  \ast\partial^{\alpha}%
\varphi_{\sqrt{\varepsilon}}\right)  (x),\\
&  =\left\langle \Lambda\left(  f\right)  ,\left\{  y\mapsto\partial^{\alpha
}\varphi_{\sqrt{\varepsilon}}\left(  x-y\right)  \right\}  \right\rangle .
\end{align*}
Consider $K$ a compact subset of $\mathbb{R}^{m}$. As $\operatorname*{supp}%
\varphi_{\sqrt{\varepsilon}}$ decrease to $\left\{  0\right\}  $ for
$\varepsilon\rightarrow0$, there exists a compact $K^{\prime}$ such that
\[
\forall x\in K,\;\;\forall\varepsilon\in\left(  0,1\right]
,\;\;\;\;\;\operatorname*{supp}\left(  \partial^{\alpha}\left(  y\mapsto
\varphi_{\sqrt{\varepsilon}}\left(  x-y\right)  \right)  \right)  \subset
K^{\prime}.
\]
The map
\[
\Theta:\mathcal{D}_{J}\left(  \mathbb{R}^{n}\right)  \times\mathcal{D}%
_{K^{\prime}}\left(  \mathbb{R}^{m}\right)  ,\;\;\;\;\left(  f,\varphi\right)
\rightarrow\left\langle \Lambda\left(  f\right)  ,\varphi\left(
x-\cdot\right)  \right\rangle
\]
is a bilinear map separately continuous since $\Lambda$ is continuous. As
$\mathcal{D}_{J}\left(  \mathbb{R}^{n}\right)  $ and $\mathcal{D}_{K^{\prime}%
}\left(  \mathbb{R}^{m}\right)  $ are Frechet spaces, $\Theta$ is globally
continuous. There exists $C>0$,\ $l_{1}\in\mathbb{N}$,\ $l_{2}\in\mathbb{N}$
such that
\[
\forall\left(  f,\varphi\right)  \in\mathcal{D}_{J}\left(  \mathbb{R}%
^{n}\right)  \times\mathcal{D}_{K^{\prime}}\left(  \mathbb{R}^{m}\right)
,\;\;\;\;\left|  \left\langle \Lambda\left(  f\right)  ,\varphi\right\rangle
\right|  \leq CP_{J,l_{1}}(f)P_{K^{\prime},l_{2}}(\varphi\left(
x-\cdot\right)  ).
\]
In particular, for $l\in\mathbb{N}$ and $\alpha\in\mathbb{N}^{m}$ with
$\left|  \alpha\right|  \leq l$, we have
\begin{equation}
\left|  \left\langle \Lambda\left(  f\right)  ,\partial^{\alpha}%
\varphi_{\varepsilon}\left(  x-\cdot\right)  \right\rangle \right|  \leq
CP_{J,l_{1}}(f)P_{K^{\prime},l_{2}}(\partial^{\alpha}\varphi_{\sqrt
{\varepsilon}}\left(  x-\cdot\right)  ), \label{GSTProofSM}%
\end{equation}
and $P_{K^{\prime},l_{2}}(\partial^{\alpha}\varphi_{\sqrt{\varepsilon}}\left(
x-\cdot\right)  )\leq P_{K^{\prime},l_{2}+l}(\partial^{\alpha}\varphi
_{\sqrt{\varepsilon}}\left(  x-\cdot\right)  )$.

Let us recall that
\[
\partial^{\alpha}\varphi_{\sqrt{\varepsilon}}\left(  x-\cdot\right)
=\partial^{\alpha}\left\{  y\mapsto\left(  \sqrt{\varepsilon}\right)
^{-m}\varphi\left(  \left(  x-y\right)  /\sqrt{\varepsilon}\right)
\kappa\left(  \left|  \ln\varepsilon\right|  \left(  x-y\right)  \right)
\right\}  .
\]
By induction on $\left|  \alpha\right|  $ and using the boundeness of
$\varphi$, $\kappa$ and their derivatives on $\mathbb{R}^{m}$, we can show
that there exists a constant $C_{1}$, depending on $\left|  \alpha\right|  $,
$\varphi$ and $\kappa$ and their derivatives but not on $\varepsilon$, such
that
\[
\sup_{y\in K^{\prime}}\left|  \partial^{\alpha}\left\{  y\mapsto\varphi
_{\sqrt{\varepsilon}}\left(  x-y\right)  \right\}  \right|  \leq C_{1}%
^{\prime}\left(  \sqrt{\varepsilon}\right)  ^{-\left(  m+\left|
\alpha\right|  +1\right)  }.
\]
It follows that there exists a constant $C_{2}$ (independent of $\varepsilon$)
such that
\[
P_{K^{\prime},l_{2}+l}(\varphi_{\varepsilon}\left(  x-\cdot\right)  )\leq
C_{2}\left(  \sqrt{\varepsilon}\right)  ^{-\left(  m+l_{2}+l+1\right)  }.
\]
Putting this result in relation (\ref{GSTProofSM}), we finally get the
existence of a constant $C_{3}$ (independent of $\varepsilon$) such that
\[
p_{_{K,l}}(L_{\varepsilon\left.  \mathcal{D}_{J}\right|  }(f))=\sup_{x\in
K,\;\left|  \alpha\right|  \leq l}\left|  \left\langle L\left(  f\right)
,\partial^{\alpha}\varphi_{\varepsilon}\left(  x-\cdot\right)  \right\rangle
\right|  \leq C_{3}\varepsilon^{-\frac{m+l_{2}+l+1}{2}}P_{J,l_{1}}(f).
\]
The sequence $r\left(  \cdot\right)  =\left\{  l\mapsto\frac{m+l_{2}+l+1}%
{2}\right\}  $ satisfies $\lim_{l\rightarrow+\infty}\left(  r(l)/l\right)
=1/2$ showing our claim.$\blacksquare$~

\subsection{Proof of proposition \ref{GSTEequaGD}}

We first have the following:

\begin{lemma}
\label{GSTLmmGDequa}For all $T\in\mathcal{D}^{\prime}\left(  \mathbb{R}%
^{m}\right)  $ $\left[  \left(  T\ast\varphi_{\sqrt{\varepsilon}}\right)
_{\varepsilon}\right]  $ is equal to $T$ in the generalized distribution sense.
\end{lemma}

\begin{proof}
Take $T\in\mathcal{D}^{\prime}\left(  \mathbb{R}^{m}\right)  $ and
$g\in\mathcal{D}\left(  \mathbb{R}^{m}\right)  $, with $K=\operatorname*{supp}%
g$. Set, for $\varepsilon\in\left(  0,1\right]  $,
\[
A_{\sqrt{\varepsilon}}=\int_{K}\left(  T\ast\varphi_{\sqrt{\varepsilon}%
}\right)  \left(  x\right)  g\left(  x\right)  \,\mathrm{d}x=\int
_{K}\left\langle T,\varphi_{\sqrt{\varepsilon}}\left(  x-\cdot\right)
\right\rangle g\left(  x\right)  \,\mathrm{d}x.
\]
As $\operatorname*{supp}\varphi_{\sqrt{\varepsilon}}$ decrease to $\left\{
0\right\}  $ for $\varepsilon\rightarrow0$, there exists a relatively compact
open subset $\Omega$ such that
\[
\forall x\in K,\;\;\forall\varepsilon\in\left(  0,1\right]
,\;\;\;\;\;\operatorname*{supp}\left(  y\mapsto\varphi_{\sqrt{\varepsilon}%
}\left(  x-y\right)  \right)  \subset\Omega.
\]
There exists $f$ continuous with compact support and $\alpha\in\mathbb{N}^{m}$
such that $T_{\left|  \Omega\right.  }=\partial^{\alpha}f$. This implies that
$\left\langle T,\varphi_{\sqrt{\varepsilon}}\left(  x-\cdot\right)
\right\rangle =\left\langle \partial^{\alpha}f,\varphi_{\sqrt{\varepsilon}%
}\left(  x-\cdot\right)  \right\rangle $ and
\[
\left(  T\ast\varphi_{\sqrt{\varepsilon}}\right)  \left(  x\right)  =\left(
\partial^{\alpha}f\ast\varphi_{\sqrt{\varepsilon}}\right)  \left(  x\right)
=\partial^{\alpha}\left(  f\ast\varphi_{\sqrt{\varepsilon}}\right)  \left(
x\right)  .
\]
By integration by part ($g$ is compactly supported) it follows that
\[
A_{\sqrt{\varepsilon}}=\int_{K}\partial^{\alpha}\left(  f\ast\varphi
_{\sqrt{\varepsilon}}\right)  \left(  x\right)  g\left(  x\right)
\,\mathrm{d}x=\left(  -1\right)  ^{\left|  \alpha\right|  }\int_{K}\left(
f\ast\varphi_{\sqrt{\varepsilon}}\right)  \left(  x\right)  \partial^{\alpha
}g\left(  x\right)  \,\mathrm{d}x.
\]

Consider now an integer $k$ and $\beta\in\mathbb{N}^{m}$ such that
$\beta=\beta_{1}+\ldots+\beta_{m}$ with $\beta_{j}\geq k$, for each
$j\in\left\{  1,\ldots,m\right\}  $. We consider $F_{\beta}$ a function such
that $\partial^{\beta}F_{\beta}=f$, which exists since $f$ is continuous. This
function is at least of class $\mathrm{C}^{k}$. We have
\begin{align*}
A_{\sqrt{\varepsilon}}  &  =\left(  -1\right)  ^{\left|  \alpha\right|  }%
\int_{K}\left(  \partial^{\beta}F_{\beta}\ast\varphi_{\sqrt{\varepsilon}%
}\right)  \left(  x\right)  \partial^{\alpha}g\left(  x\right)  \,\mathrm{d}%
x,\\
&  =\left(  -1\right)  ^{\left|  \alpha\right|  +\left|  \beta\right|  }%
\int_{K}\left(  F_{\beta}\ast\varphi_{\sqrt{\varepsilon}}\right)  \left(
x\right)  \partial^{\alpha+\beta}g\left(  x\right)  \,\mathrm{d}x,\\
\left\langle T,g\right\rangle  &  =\left\langle \partial^{\alpha
}f,g\right\rangle =\left\langle \partial^{\alpha+\beta}F_{\beta}%
,g\right\rangle =\left(  -1\right)  ^{\left|  \alpha\right|  +\left|
\beta\right|  }\left\langle F_{\beta},\partial^{\alpha+\beta}g\right\rangle
,\\
&  =\left(  -1\right)  ^{\left|  \alpha\right|  +\left|  \beta\right|  }%
\int_{K}\left(  F_{\beta}\right)  \left(  x\right)  \partial^{\alpha+\beta
}g\left(  x\right)  \,\mathrm{d}x.
\end{align*}
Then
\[
\left\langle T\ast\varphi_{\sqrt{\varepsilon}},g\right\rangle -\left\langle
T,g\right\rangle =\left(  -1\right)  ^{\left|  \alpha\right|  +\left|
\beta\right|  }\int_{K}\left(  \left(  F_{\beta}\ast\varphi_{\sqrt
{\varepsilon}}\right)  \left(  x\right)  -\left(  F_{\beta}\right)  \left(
x\right)  \right)  \partial^{\alpha+\beta}g\left(  x\right)  \,\mathrm{d}x.
\]
An adaptation (and simplification) of the proof of lemma \ref{LmnMTh4} shows
that
\[
\left(  F_{\beta}\ast\varphi_{\sqrt{\varepsilon}}\right)  \left(  x\right)
-\left(  F_{\beta}\right)  \left(  x\right)  =\mathrm{O}\left(  \sqrt
{\varepsilon}^{k}\right)  \text{ for }\varepsilon\rightarrow0\text{.}%
\]
As $g$ is compactly supported, this last relation leads to
\[
\left\langle T\ast\varphi_{\sqrt{\varepsilon}},g\right\rangle -\left\langle
T,g\right\rangle =\mathrm{O}\left(  \sqrt{\varepsilon}^{k}\right)  \text{ for
}\varepsilon\rightarrow0.
\]
Since $k$ is arbitrary, our claim follows.
\end{proof}

This lemma implies that for all $f\in\mathcal{D}\left(  \mathbb{R}^{n}\right)
$, $\left[  \left(  L_{\varepsilon}(f)\right)  _{\varepsilon}\right]  =\left[
\left(  \Lambda\left(  f\right)  \ast\varphi_{\sqrt{\varepsilon}}\right)
_{\varepsilon}\right]  $ is equal to $\Lambda\left(  f\right)  $ in the
generalized distribution sense. On the other hand, according to theorem
\ref{ThmGST}, $\left[  \left(  L_{\varepsilon}(f)\right)  _{\varepsilon
}\right]  =\widetilde{H}_{L}\left(  f\right)  $ where $\widetilde{H}_{L}$ is
the integral operator associated to the canonical extension of $\left(
L_{\varepsilon}\right)  _{\varepsilon}$.\ This ends the proof.

\end{document}